\documentclass{article}
\usepackage{nameref}
\usepackage{float}
\usepackage{pdflscape}
\usepackage[plainpages=false,pdfpagelabels]{hyperref}
\usepackage[english]{babel}
\usepackage[utf8]{inputenc}
\usepackage{latexsym}
\usepackage{textcomp}
\usepackage{amsmath,amsthm}
\usepackage{amsfonts}
\usepackage{epsfig}
\usepackage{amssymb}
\usepackage{longtable}
\usepackage{tabularx,colortbl}
\usepackage[table]{xcolor}
\usepackage{graphicx,epsfig,psfrag}
\usepackage{color}
\usepackage{eso-pic}
\usepackage{pdflscape}
\usepackage{dcolumn}
\usepackage{flafter}
\usepackage[mathscr]{euscript}
\usepackage{array, lscape}
\usepackage{pb-diagram}
\usepackage{setspace}
\usepackage{apacite}
\usepackage{rotating}
\usepackage[all]{hypcap}
\usepackage[font=footnotesize]{caption}

\newtheorem{corollary}{Corollary}

\newtheorem{lemma}{Lemma}

\newtheorem{thm}{Theorem}
\begin{document}
\newcommand*\samethanks[1][\value{footnote}]{\footnotemark[#1]}
% **** --------------------------------------------------------------------------------
\title{A derivation of the optimal answer-copying index and some applications\thanks{Corresponding author: Mauricio Romero. e-mail: \href{mailto:mtromero@ucsd.edu}{mtromero@ucsd.edu}. The authors would like to thank the ICFES for financial support, three anonymous referees, Nicola Persico, Decio Coviello and Julian Mari\~no and his group of statisticians for valuable comments and suggestions.}}
\author{Mauricio Romero\thanks{University of California - San Diego.} \thanks{Quantil $\mid$ Matem\'{a}ticas Aplicadas.} \and \'Alvaro Riascos \thanks{Universidad de los Andes.} \samethanks[3]   \and Diego Jara\samethanks[3]}% ////
\maketitle
\begin{abstract}
\footnotesize{Multiple-choice exams are frequently used as an efficient and objective method to assess learning but they are more vulnerable to answer-copying than tests based on open questions. Several statistical tests (known as indices in the literature) have been proposed to detect cheating; however, to the best of our knowledge they all lack mathematical support that guarantees optimality in any sense. We partially fill this void by deriving the uniform most powerful (UMP) under the assumption that the response distribution is known. In practice, however, we must estimate a behavioral model that yields a response distribution for each question. We calculate the empirical type-I and type-II error rates for several indices that assume different behavioral models using simulations based on real data from twelve nationwide multiple-choice exams taken by 5th and 9th graders in Colombia. We find that the index with the highest power among those studied, subject to the restriction of preserving the type-I error, is one based on the work of \citeA{Wollack1997} and \citeA{Linden2006} and is superior to the indices studied and developed by  \citeA{Wesolowsky2000} and \citeA{Frary1977}.  We compare the results of applying this index to all 12 exams and find that examination rooms with stricter proctoring have a lower level of copying. Finally, a Bonferroni correction to control for the false positive rate is proposed to detect massive cheating.}

\textbf{Key Words:} $\omega$ Index , Answer Copying, False Discovery Rate, Neyman-Pearson's Lemma.

\textbf{JEL Clasification:} C19, I20
\end{abstract}

\section{Introduction}

Multiple-choice exams are frequently used as an efficient and objective way of evaluating knowledge. Nevertheless, they are more vulnerable to answer copying than tests based on open questions. Answer-copy indices provide a statistical tool for detecting cheating by examining suspiciously similar response patterns between two students. However, these indices have three problems. First, similar answer patterns between a pair of students could be justified without answer copying. For example, two individuals with very similar educational background are likely to provide similar answers. The second problem is that a statistical test (an index) is by no means a conclusive basis for accusing someone of copying, since it is impossible to completely eliminate type-I errors. In other words, it is possible that two individuals share the same response pattern by chance. Finally, every index assumes responses are stochastic. If the assumed probability distribution is incorrect, the index can lead to incorrect conclusions. Furthermore, all the indices in the literature are ad-hoc and there are no theoretical results that support the use of one index over the other.

%Generally speaking, these tests should only be used as a first alert, and not as those sole or final proof of cheating. 
%

\citeA{Wollack2003} compares several indices and finds that among those that preserve size the $\omega$ index is the most powerful one. However, the set of indices studied is not comprehensive and in particular does not include the index developed by \citeA{Wesolowsky2000}.

%The $\kappa$ index, however, is not among those compared by \citeA{Wollack2003}, thus making it hard for the ICFES to justify a migration from one index to the other. Additionally, 

Thus there are two gaps in the literature that this article seeks to fill. First, it provides theoretical foundations that validate the use of indices that reject the null hypothesis of no cheating for a large number of identical answers under the assumption that student responses are stochastic. 

Second, it compares the type-I and type-II error rates of the $\omega$ and $\gamma$ indices for answer copy detection, based on the work of \citeA{Wollack1997}\footnote{In this article we use a version closer to the work of \citeA{Linden2006}.} and \citeA{Wesolowsky2000} respectively\footnote{Both indices are refinements of the indices first developed by \citeA{Frary1977}.}. Using Monte Carlo simulations and data from the SABER tests taken by 5th and 9th graders in Colombia in May and October of 2009 we find that the conditional version of the standardized index first developed by \citeA{Wollack1997} is the most powerful among those that respect size. 

We compare the results of applying the index to examination rooms with different strategies to control cheating. We find a negative correlation between the level of proctoring and the prevalence of copying. We also find a lower prevalence of copying in examination rooms where students answer different portions of the test at the same time compared to examination rooms where all students answer the same portion of the test at the same time. These results have at least two possible interpretations: they could be interpreted as evidence that the index is indeed detecting cheating, or, alternatively, if one believes that the index can be used as a reliable measure of cheating, these results can be interpreted as estimates of how effective current strategies for cheating-prevention are.  However, the results of these two exercises must be taken cautiously as they are not the result of a randomized experiment and therefore might be biased due to unobservable factors.

Our article has a fourth contribution. We outline a procedure for detecting massive cheating. These indices detect individual cheating, but do not constitute definitive proof of copying given their statistical nature. They are merely intended to raise flags. In Colombia an index is used to search for examination rooms with a large number of flags (i.e. a large proportion of students guilty of copying according to the index). When a large number of flags are raised, every student in the examination room must retake the test under stricter surveillance conditions. The appropriate way to search for a large number of flags is to test multiple hypothesis at the same time, but these procedures often results in low statistical power. We apply a Bonferroni correction outlined by \citeA{Benjamini1995} to detect multiple cheating while controlling for the false positive rate.  The application is straightforward and we think this is a useful tool for flagging possible examination rooms where massive cheating might have occurred. This information could be used, as in Colombia, to make entire examination rooms retake an exam under stricter surveillance conditions.

The article is organized as follows. The second section derives an optimal statistical test (index) to detect answer copying using the Neyman-Pearson's Lemma. The third section presents two of the most widely used indices, which are based on the work of \citeA{Wollack1997}, \citeA{Frary1977}, \citeA{Wesolowsky2000}, and \citeA{Sotaridona2006}. The fourth section presents a brief summary of the data used and is followed by a section that presents the methodology of the Monte Carlo simulations used to find the empirical type-I and type-II error rates (to test which behavioral model gives the best results) and its results.  Section six analyses the correlation between different strategies to control cheating and the prevalence of cheating according to the index and section seven presents the results of using a Bonferroni correction to calculate the prevalence of massive cheating. Finally the last section concludes.

\section{Applying Neyman-Pearson's to answer copying}

%All answer copying indices try to find if the similarities in two answer patterns are suspiciously unlikely. In other words 

It is normal for two answer patterns to have similarities by chance. Answer-copying indices try to detect similarities that are so unlikely  to happen naturally that answer-copying becomes a more natural explanation than chance. Most answer-copy indices are calculated by counting the number of identical answers between the test taker suspected of copying and the test taker suspected of providing answers\footnote{For examples see \citeA{Linden2004,Linden2006,Sotaridona2003,Sotaridona2002,Sotaridona2006,Holland1996,Frary1977,Cohen1960,Bellezza1989,Angoff1974, Wesolowsky2000,Wollack1997}}.  In all these indices the null hypothesis is the same: there is no cheating.

All these indices are ad-hoc since they are not derived to be optimal in any sense. To the extent of the authors' knowledge, this article presents the first effort to rationalize the use of these indices to detect answer copying using the Neyman-Pearson's Lemma (NPL) \cite{Neyman1933} resulting in the uniformly most powerful (UMP) test (index), assuming we know the underlying probability of responses for each individual in each question. However, we must turn to empirical data to find the performance of each index since different behavioral models result in different response distributions.

First, let us state the problem formally. Let us assume that there are $N$ questions and $n$ alternatives for each question. We are interested in testing whether the individual who cheated (denoted by $c$) copied from the individual who supposedly provided the answers (denoted by $s$). Let $\gamma_{cs}$ be the number of questions that $c$ copied from $s$. The objective is to test the following hypotheses:
\begin{align*}
	H_0: \gamma_{cs}=0\\
	H_1: \gamma_{cs}>0\\
\end{align*} 
%The results that will follow do not depend on how we model peoples response behavior. , we only need $\pi_{iv}^j$ for all $i$,$v$ and $j$. 

Let $I_{csi}$ be equal to one when individuals $c$ and $s$ have the same answer to question $i$ and zero otherwise. Then, the number of common answers between $c$ and $s$ can be expressed as:
\begin{equation}
M_{cs}=\sum_{i=1}^N I_{csi} .
\end{equation}

Under the null hypothesis $M_{cs}$ is the sum of $N$ independent Bernoulli random variables, each with a different probability of success $\pi_i$, equal to the probability that individual $c$ has the same answer as individual $s$ in question $i$. The distribution of $M_{cs}$ is known as a poisson binomial distribution. Let $B(\pi_1,...,\pi_N)$ be such distribution and $f_N(x;\pi_1,...,\pi_N)$ the probability mass function (pmf) at $x$. Notice that if $\pi_1=\pi_2=...=\pi_N=\pi$ then the poisson binomial distribution reduces to a standard binomial distribution.

Now, let $A$ denote the set of questions that student $c$ copied from $s$. Then if $|A|=k$, it means that $\gamma_{cs}=k$, and $M_{cs}$ has the following probability mass function (pmf) $\hat{f}_N(x;\pi_1,...,\pi_N, A)$, where we define $ \hat{f}_N(x;\pi_1,...,\pi_N, A) \doteq f_N(x,\pi_{1}',..,\pi_{N}')$ such that 
$$\pi_i'=\begin{cases}
1 & \text{ if } i \in A \\
\pi_i & \text{ if } i \not\in A \\
\end{cases}$$

%$$f(x;A)=\begin{cases}
%0 & \text{ if } x<k \text{ or } x>N \\
%\hat{f}(x,\pi_1,...,\pi_N; A) & \text{ if } k \leq x \leq N \\
%\end{cases}$$

For example, say that there are 50 questions and that the students copied questions 1, 10 and 50, i.e. $A=\lbrace 1, 10, 50 \rbrace$ then $$\hat{f}_N(x;\pi_1,...,\pi_N, A)=f_N(x;1,\pi_2,...,\pi_9,1,\pi_{11},...,\pi_{49},1).$$

Before we continue let us state Neyman-Pearson's Lemma (NPL):

\begin{thm}{Neyman-Pearson's Lemma \cite{casella2002}}

Consider testing $H_0: \theta=\theta_0$ against $H_1:\theta=\theta_1$ where the pmf is $f(\textbf{x}|\theta_i)$, $i=0,1$, using a statistical test (index) with rejection region $R$ that satisfies 

\begin{equation}\label{eq1}
\begin{aligned}
\textbf{x} \in R & \text{ if } f(x|\theta_1) > f(x|\theta_0) k\\
\textbf{x} \in R^c & \text{ if } f(x|\theta_1) < f(x|\theta_0) k
\end{aligned}
\end{equation}

for some $k \geq 0$, and

\begin{equation}\label{eq2}
\alpha=P_{H_0}(\textbf{X} \in R)
\end{equation}

Then 

\begin{enumerate}
\item (Sufficiency) Any test (index) that satisfies equations \ref{eq1} and \ref{eq2} is a UMP level $\alpha$ test (index).
\item (Necessity) If there exists a test (index) satisfying equations \ref{eq1} and \ref{eq2} with $k>0$, then
every UMP level $\alpha$ test (index) is a size $\alpha$ test (index) - satisfies \ref{eq2} - and every UMP level $\alpha$ test (index) satisfies \ref{eq1} except perhaps on a set $A$ such that $P_{H_0}(\textbf{X} \in A) = P_{H_1}(\textbf{X} \in A)=0$.
\end{enumerate}

the test (index) is the uniformly most powerful (UMP) level $\alpha$ test (index).

%   Then the test with the highest power among all tests with significance level $\alpha$ has rejection region
%
%
%
%$$R(\textbf{x}) = \lbrace \textbf{x} : \lambda(\textbf{x}) \geq k \rbrace$$
%
%where
%
%$$\lambda(\textbf{x})=\frac{f(\textbf{x}|\theta_1)}{f(\textbf{x}|\theta_0)}$$
%
%is the likelihood ratio and $k$ is a constant such that
%
%$$P_{H_0}\left( \lambda(\textbf{x}) \leq k \right)=\alpha$$
%

\end{thm} 

In this context, let us apply the NPL to the simple hypothesis test $H_0: A=A_0$ and $H_1: A=A_1$, where $A_0=\emptyset$ (i.e. there is no cheating) and $A_1$ is a set of questions, to get the UMP test. If in the data we observe $x$ questions answered equally by individuals $c$ and $s$ then the likelihood ratio test\footnote{Notice that NPL implies that a likelihood ratio test is the uniformly most powerful (UMP) test for simple hypothesis testing.} would be:

$$\lambda^A(x)=\frac{\hat{f}_N(x;\pi_1,...,\pi_N, A)}{f_N(x;\pi_1,...,\pi_N)}$$

Now we must find the critical value of the test. In other words, we need the greatest value $c$ such that under the null we have:

$$1-P_{H_0} \left( \frac{\hat{f}_N(x;\pi_1,...,\pi_N, A)}{f(x;\pi_1,...,\pi_N)}<c \right)= P_{H_0} \left(\frac{\hat{f}_N(x;\pi_1,...,\pi_N, A)}{f_N(x;\pi_1,...,\pi_N)}>c \right)\leq \alpha$$

For any given pair of simple hypotheses ($H_0: A=A_0$ , $H_1: A=A_1$) we know how to find the UMP (by using the NPL). The following lemma will allow us to find the UMP for more complex alternative hypothesis (e.g. $H_1: \lbrace A : |A| \geq 1 \rbrace $).

\begin{lemma}
$\lambda^A(x)=\frac{\hat{f}_N(x;\pi_1,...,\pi_N, A)}{f_N(x,\pi_1,...,\pi_N)}$ is increasing in $x\in \lbrace 0,...,N\rbrace$ for all $A$. 
\end{lemma}
Before we present the proof we must first recall some useful results proved by \citeA{Wang1993}.
\begin{thm}[Theorem 2 in \citeA{Wang1993}]
The pmf of a poisson binomial satisfies the following inequality:
$$f_{N}(x;\pi_1,\pi_2,...,\pi_N)^2 > C(x) f_{N}(x+1;\pi_1,\pi_2,...,\pi_N)f_{N}(x-1;\pi_1,\pi_2,...,\pi_N)$$
where $C(x)=\max \left(\frac{x+1}{x},\frac{N-x+1}{N-x} \right)$
\end{thm}
which has as an immediate corollary
\begin{corollary}\label{color1}

The pmf of a poisson binomial satisfies the following inequality:
$$f_{N}(x;\pi_1,\pi_2,...,\pi_N)^2 \geq f_{N}(x+1;\pi_1,\pi_2,...,\pi_N)f_{N}(x-1;\pi_1,\pi_2,...,\pi_N) $$
\end{corollary}

Now we are ready to prove the lemma:

\begin{proof}[Proof of Lemma 1]

We consider the case $|A|=1$, given that the proof for the case when $|A|>1$ can be obtained by induction. Without loss of generality, assume $A=\{1\}$. The numerator in the lemma's quotient is $0$ for $x=0$, so we proceed to prove monotonicity $\lambda^A(x)$ in $x$ for $x\geq1$. Likewise, the case $N=1$ follows trivially, so we assume $N>1$.

For simplicity, we call $g(x)=f_{N-1}(x;\pi_2,\ldots,\pi_N)$. First, note that 
\begin{equation*}
\hat{f}_N(x;\pi_1,\ldots,\pi_N;A)=g(x-1).
\end{equation*}
Second, corollary \ref{color1} states that $g(x-1)g(x+1) \leq g(x)^2$. Third, we can write $f_N(x;\pi_2,\ldots,\pi_N)=\pi_1g(x-1)+(1-\pi_1)g(x)$. With these observations we have
\begin{align*}
\frac{\hat{f}_N(x;\pi_1,\ldots,\pi_N;A)}{f_N(x;\pi_1,\ldots,\pi_N)} &= \frac{g(x-1)}{\pi_1g(x-1)+(1-\pi_1)g(x)}\times\frac{\pi_1g(x)+(1-\pi_1)g(x+1)}{\pi_1g(x)+(1-\pi_1)g(x+1)}\\
&\leq \frac{\pi_1g(x)g(x-1)+(1-\pi_1)g(x)^2}{[\pi_1g(x-1)+(1-\pi_1)g(x)][\pi_1g(x)+(1-\pi_1)g(x+1)]}\\
&= \frac{g(x)}{\pi_1g(x)+(1-\pi_1)g(x+1)}\\
&= \frac{\hat{f}_N(x+1;\pi_1,\ldots,\pi_N;A)}{f_N(x+1;\pi_1,\ldots,\pi_N)}.
\end{align*}

\end{proof}

\bigskip

Given that $\frac{\hat{f}_N(x;\pi_1,...,\pi_N; A)}{f_N(x;\pi_1,...,\pi_N)}$ is increasing in $x$ for all $A$ then we have that for every $c$ there exists a $k^*$ such that $P_{H_0} \left( \frac{\hat{f}_N(x;\pi_1,...,\pi_N, A)}{f_N(x;\pi_1,...,\pi_N)}<c \right)= \sum_{w=0}^{k^*} f_N(w,\pi_1,...,\pi_N)$. In particular for a given level $\alpha$ of the test we can find $k^*$ such that 
$$1-P_{H_0} \left( \frac{\hat{f}_N(x;\pi_1,...,\pi_N, A)}{f(x;\pi_1,...,\pi_N)}<c \right)=\sum_{w=0}^{k^*} f(w,\pi_1,...,\pi_N) \leq \alpha$$

Then, if we reject the null hypothesis when $M_{cs}>k^*$, we get the UMP for a particular set A. However, the rejection region is the same for all $A$, thus if we reject the null hypothesis when $M_{cs}>k^*$, we get the UMP for all $A$ such that $|A|\geq 1$. This justifies the use of indices that reject the null hypothesis for large values of $M_{cs}$. 

However, $\pi_{i}$ must be estimated somehow (it was taken as known in this section and thus in empirical applications we do not have the UMP) and different methods yield different results. We now turn to the data to find out which index performs better in practice. 

\citeA{Frary1977} in a seminal article developed the first indices, known as $g_1$ and $g_2$, that reject the null hypothesis for large values of $M_{cs}$. \citeA{Wollack1997}, \citeA{Linden2006} and \citeA{Wesolowsky2000} have proposed further refinements of \citeA{Frary1977} methods. We will evaluate the performance of these indices in practice. 

% and thus it seems natural to see how this index performs in practice. However, the $\omega$ index developed by \citeA{Wollack1997} performs better than the $g_2$ \cite{Wollack2006} so we turn our attention to it and its non-standarized version developed \citeA{Linden2006}. We also focus on two other variants of the $g_2$ developed by \citeA{Wesolowsky2000}. 
 %Furthermore, lets assume that each student $j$ has a probability $\pi_{iv}^j$ of answering option $v$ on question $i$.
 
\section{Copy Indices}

Let us assume that student $j$ has a probability $\pi_{iv}^j$ of answering option $v$ on question $i$. The probability that two students have the same answer on question $i$ ($\pi_i$) can be calculated in two ways. First, assuming independent answers, the probability of obtaining the same answer is $\pi_i=\sum_{v=1}^n \pi_{iv}^c \pi_{iv}^s$. 

Second, we could think of the answers of individual $s$ as being fixed, as if he were the source of the answers and $c$ the student who copies. In the absence of cheating, conditional on the answers of $s$, the probability that individual $c$ has the same answer as individual $s$ in question $i$ is $\pi_i= \pi_{iv_s}^c$, where $\pi_{iv_s}^c$ is the probability that individual $c$ answered option $v_s$ which was chosen by $s$ in question $i$.

A discussion of these two approaches is given in \citeA{Frary1977} and \citeA{Linden2006}. The first is known as the unconditional index and is symmetric in the sense that the choice of who is $s$ and who is $c$ is irrelevant since $\pi_i$ is the same either way. The second is known as the conditional index and it is not symmetric opening the possibility that the index rejects the null hypothesis that student $a$ copied from student $b$ but not rejecting the null hypothesis that $b$ copied from $a$. The details of each situation determine which approach is appropriate. If we believe students copied from each other or answered the test jointly then a conditional index is undesirable, but if we believe that a student is the source (for whatever reason) of answers but did not collaborate with the cheater, then a conditional index might be more appropriate. We study both conditional and unconditional indices.

Indices vary along three dimensions. The first dimension is how they estimate $\pi_{iv}^j$. The second is whether they are a conditional or an unconditional index. Finally, they vary how critical values are calculated. They either use the exact distribution (a poisson binomial distribution) or a normal distribution, by applying some version of the central limit theorem. 

In order to use the central limit theorem in this context recall $M_{cs}$ is the sum of $N$ Bernoulli variables and has mean $\sum_{i=1}^N \pi_i$ and variance $\sum_{i=1}^N \pi_i(1-\pi_i)$. Thus $\frac{M_{cs}-\sum_{i=1}^N \pi_i}{\sqrt{\sum_{i=1}^N \pi_i(1-\pi_i)}}$ converges in distribution to a standard normal distribution as $N$ goes to infinity. There are two advantages to the normal approximation. First critical values are easier to calculate and more precise (computationally) and second it allows for a finer choice of critical values.

As mentioned before, \citeA{Frary1977} developed the first indices, known as $g_1$ and $g_2$, that reject the null hypothesis for large values of $M_{cs}$. However, both \citeA{Wesolowsky2000} and \citeA{Wollack2003} show that variations of the original method proposed by \citeA{Frary1977} yield superior results, and in this article we study the indices they developed. The first variation is the $\omega$ index developed by \citeA{Wollack1997} that assumes there is an underlying nominal response model. The second variation is the $\gamma$ index developed by \citeA{Wesolowsky2000} based on a variation of \citeA{Frary1977} work.

\subsection{$\omega$ index}
The $\omega$ index is based on the work of \citeA{Wollack1997} and assumes a nominal response model that allows the probability of answering a given option to vary across questions and individuals. As before, let $N$ be the number of questions and $n$ the number of alternatives for answering each question. Suppose that an individual with skill $\theta_j$, who does not copy, responds with probability $\pi_{iv}$ for option $v$ to question $i$. In other words:

% We assume a nominal response model, since this is the most widely used polytomous response model \cite{Thissen1984}. \citeA{DeAyalaWinter1989} shows that the nominal response model is superior to the three-parameter logistic model, in the sense that it estimates ability more accurately. 
\begin{equation}\label{resp_nom}
\pi_{iv}(\theta_j)=\frac{e^{\xi_{iv}+\lambda_{iv}\theta_j}}{\sum_{h=1}^m e^{\xi_{ih}+\lambda_{ih}\theta_j}},
\end{equation}

where $\xi_{iv}$ y $\lambda_{iv}$ are model parameters and are known as the intercept and slope, respectively. The intercept and slope can vary across questions. The parameters of the questions ($\xi_{iv}$ and $\lambda_{iv}$) are estimated using marginal maximum likelihood, while ability is estimated using the EAP method (Expected A Posteriori). The estimation is performed using the \textit{rirt} package in R \cite{Germain}\footnote{The package \textit {rirt} can be found on: \url{http://libirt.sourceforge.net/}.}. It is necessary to estimate ability as the proportion of correct answers taking into account that a correct answer to a ``difficult'' question indicates a higher ability than a correct answer to a ``simple'' question.  More information on marginal maximum likelihood and EAP can be found in \citeA{Linden1997} and \citeA{Hambleton1991}.

Let $\omega_1$ and $\omega_2$ be the unconditional and conditional (exact) versions of this index (following somewhat the $g_1$ and $g_2$ notation of \citeA{Frary1977}) and let $\omega_1^s$ and $\omega_2^s$ be the standardized versions (i.e. they use the normal distribution to find the critical values of the index). 

\subsection{$\gamma$ index}

The indices developed by \citeA{Wesolowsky2000}, which are an extension and improvement of the indices developed by \citeA{Frary1977}, assume that the probability that student $j$ has the correct answer in question $i$ is given by:

$$p_{i}=(1-(1-r_i)^{a_j})^{1/a_j},$$

where $r_i$ is the proportion of students that had the right answer in question $i$. The parameter $a_j$ is estimated by solving the equations
$$\frac{\sum_{i=1}^N p_i}{n}=c_j,$$ 

where $c_j$ is the proportion of questions answered correctly by individual $j$. Finally, we need the probability that student $j$ chooses option $v$ among those that are incorrect which is estimated as the proportion of students with an incorrect answer that chose each incorrect option. Lets denote $\gamma_1$ and $\gamma_2$ the unconditional and conditional version of this index and by $\gamma_1^s$ and $\gamma_2^s$ their standardized version respectively.

Before we compare how the different versions of the $\omega$ and the $\gamma$ index fare in practice, the following section presents the data that will be used.

\section{Data}

In Colombia, every student enrolled in 5th, 9th or 11th grade, whether attending a private or a public school, has to take a standardized multiple-choice test known as the SABER test\footnote{The tests in the 5th and 9th grade have been somewhat irregular and with students being tested every 2 to 3 years.}. These exams are intended to measure the performance of students and schools across several areas. The Instituto Colombiano para la Evaluaci\'{o}n de la Educaci\'{o}n (ICFES) is in charge of developing, distributing and applying these exams. The score of the 11th grade test is used by most universities in Colombia as an admission criterion. The ICFES also evaluates all university students during their senior year. We analyze the 5th and 9th grade tests for 2009. In total, we have 12 different exams depending on the subject, the date of the exam and the grade of the student\footnote{Each grade (5th and 9th) presents three tests: Science, Mathematics and Language. Schools that finish the academic year in December present the exam in September and schools that finish their academic year in June present the exam in May. In total there are two dates, two grades and three subjects, for a total of 12 exams.}.The following abbreviations, used by the ICFES, are used: per grade, 50 for 5th and 90 for 9th. Per area, 041 for mathematics, 042 for language and 043 for science. Per date, F1 for May and F2 for October.  For example, exam PBA9041F2 is taken by 9th graders for mathematics in October. A brief overview of each test is presented in Table \ref{tb:ov}.

%Between 2010 and 2011 the ICFES used the $\kappa$ index, based on the work of \citeA{Sotaridona2006}, as a mechanism for detecting cheating in these exams. Using data from the ICFES exams \citeA{Jara2010} showed that in practice the $\kappa$ index has a larger type-I error rate than predicted by theory suggesting that the assumptions used to construct the index are not met and that the ICFES should use a different index.  Many institutions in charge of standardized testing regularly face a similar problem and need theoretical and empirical evidence to justify the use of an index.

For each exam the database contains the answer to each question for each individual, as well as the examination room where the exam was taken. The correct answers for each exam are also available.

\begin{table}[H]
\caption{Summary statistics}
\label{tb:ov}
\begin{center}
\begin{tabularx}{1.05\textwidth}{|l|l|l|l|X|X|X|}
  \hline
Test & Subject & Grade & Month & Questions & Students & Examination Rooms \\
  \hline
     PBA5041F1  & Math & 5th   & May   &  48   &  60,099  &  3,421  \\
     PBA5041F2  & Math & 5th   & Oct &  48   &  403,624  &  31,827  \\
     PBA5042F1  & Language & 5th   & May   &  36   &  60,455  &  3,441  \\
     PBA5042F2  & Language & 5th   & Oct &  36   &  402,508  &  31,642  \\
     PBA5043F1  & Science & 5th   & May   &  48   &  60,404  &  3,432  \\
     PBA5043F2  & Science & 5th   & Oct &  48   &  405,537  &  31,833  \\
     PBA9041F1  & Math & 9th   & May   &  54   &  44,577  &  1,110  \\
     PBA9041F2  & Math & 9th   & Oct &  54   &  303,233  &  9,059  \\
     PBA9042F1  & Language & 9th   & May   &  54   &  44,876  &  1,110  \\
     PBA9042F2  & Language & 9th   & Oct &  54   &  302,781  &  9,044  \\
     PBA9043F1  & Science & 9th   & May   &  54   &  44,820  &  1,107  \\
     PBA9043F2  & Science & 9th   & Oct &  54   &  30,3723  &  9,053  \\
       \hline
\multicolumn{7}{p{1\textwidth}}{Source: ICFES. Calculations: Authors. }\\
%\multicolumn{7}{p{1\textwidth}}{Note: The number of schools corresponds to the number of examination rooms.}\\
\end{tabularx}
\end{center}
\end{table}

\section{Index Comparison}
In this section we compare the different versions of the $\omega$ and the $\gamma$ indeces. In order to do this we evaluate the type-I and type-II error rates by creating synthetic samples in which we control the level of cheating between individuals. 

\subsection{Methodology}

To find the empirical type-I error rate, individuals who could not have possibly copied from one another are paired together and tested for cheating using a particular index. This is done by pairing individuals that answered the exam in different examination rooms, thus eliminating the possibility of answer copying. The empirical type-I error rate is calculated as the proportion of pairs for which the index rejects the null hypothesis. To find the empirical type-II error rate, we take these answer-copy free pairs and simulate copying by forcing specific answers to be the same. The proportion of pairs for which the index rejects the null hypothesis is the power of the index\footnote{Recall that the power of the test is the complement of the type-II error rate, i.e. $Power=100\%-Type II Error$.}.

To make things clearer, let $c$ denote the test taker suspected of cheating, $s$ the test taker believed to have served as the source of answers. The steps taken to find the type-I error rate and the power of each index are as follows:

\begin{enumerate}
\item 100,000 pairs are picked in such a way that for each couple the individuals performed the exam in different examination rooms.
\item The answer-copy methodology is applied to these pairs, and the proportion of pairs  for which the index rejects the null hypothesis is the empirical type-I error rate estimator.
\item To calculate the power of the index, the answer pattern for individual $c$ is changed by replacing $k$ of his answers to match to those of individual $s$\footnote{For example, let us assume the answer pattern for $s$ is $ACBCDADCDAB$, which means that there were 11 questions and that he/she answered A for the first questions, C for the second questions, and so on. Also assume that the original answer pattern of $c$ without copying is $DCABCDAABCB$. Let $k$ be 5, this means and let us assume that the randomly selected questions were 1,4,5,10,11. This means that the modified (with copying) answer patterns for $c$ will be $ACACDDAABAB$.}.
\begin{enumerate}
\item The level of copy, $k$, is set, and is defined as the number of answers transferred from $s$ to $c$.
\item $k$ questions are selected randomly.
\item Individual c's answers for the $k$ questions are changed to replicate exactly those of individual $s$. Answers that were originally identical count as part of the $k$ questions being changed.
\end{enumerate}
\item We apply the answer-copy methodology to the pairs whose exams have been altered. The proportion of pairs accused of cheating is the power of the index for a copying level of $k$.
\end{enumerate}

\subsection{Results}
%All the results, including the estimated coefficients for the nominal response model, are available upon request. The results for the grid of abilities can be found at \url{www.quantil.com.co/appendixOmega_blind.pdf}.
%
%\subsection{The results without conditioning for ability}
Throughout the analysis a confidence level ($\alpha$) of 99.9\% is used and the power of the index is calculated at copying levels ($k$) of: $1,5,10,15,20,...,N$, where $N$ is the number of questions in the exam.

\subsubsection{Type-I error rate}\label{tipoI}

As can be seen in Tables \ref{tab:errotipo1general1} and \ref{tab:errotipo1general2} the $\gamma_2$, $\gamma_2^s$, and the $\omega_2$ indices have an empirical type-I error rate that is consistently above the theoretical type-I error rate of one in a thousand. The $\gamma_1$ index (developed by \citeA{Wesolowsky2000}) empirical error rate is above the theoretical one in several cases. 

Based on these results, we discard the $\gamma_2$, $\gamma_2^s$ and the $\omega_2$ indices and restrict the search for the most powerful index among $\gamma_1$, $\gamma_1^s$, $\omega_1$ and $\omega_2^s$. 
 
\begin{table}[H]
  \centering
  \caption{Type-I error rate for the $\gamma$ indices}
   \begin{tabularx}{\textwidth}{|c|c|c|c|c|c|c|c|}
    \hline
        Exam  & Subject & Grade & Month & $\gamma_1$ & $\gamma_2$ & $\gamma_1^s$ & $\gamma_2^s$  \\
        \hline
        PBA5041F1  & Mathematics & 5th   & May   &  0.66 & 2.20 & 0.41 & 0.76 \\ 
     PBA5041F2  & Mathematics & 5th   & October &  0.87 & 2.44 & 0.59 & 1.11 \\
     PBA5042F1  & Language & 5th   & May   &   1.20 & 2.18 & 0.77 & 1.16 \\ 
     PBA5042F2  & Language & 5th   & October &  1.21 & 2.36 & 0.92 & 1.49 \\
     PBA5043F1  & Science & 5th   & May   &  1.05 & 2.59 & 0.73 & 1.38 \\ 
     PBA5043F2  & Science & 5th   & October &   0.74 & 1.81 & 0.61 & 1.24 \\ 
     PBA9041F1  & Mathematics & 9th   & May   &  1.38 & 1.97 & 0.96 & 1.26 \\
     PBA9041F2  & Mathematics & 9th   & October & 2.15 & 2.14 & 1.69 & 1.53 \\
     PBA9042F1  & Language & 9th   & May   &  0.85 & 2.24 & 0.56 & 1.04 \\ 
     PBA9042F2  & Language & 9th   & October & 0.84 & 1.92 & 0.59 & 1.34 \\ 
     PBA9043F1  & Science & 9th   & May   & 1.32 & 2.06 & 0.93 & 1.42 \\
     PBA9043F2  & Science & 9th   & October &  1.02 & 1.70 & 0.74 & 1.37 \\
  
        \hline
        \multicolumn{8}{p{0.9\textwidth}}{Source: ICFES. Calculations: Authors.}  \\      
        \multicolumn{8}{p{0.9\textwidth}}{Number of copy-free couples accused of copying (for every 1,000 pairs) at a 99.9\% confidence level}\\
    \end{tabularx}
  \label{tab:errotipo1general1}
\end{table}
\begin{table}[H]
  \centering
  \caption{Type-I error rate for the $\omega$ indices}
    \begin{tabularx}{\textwidth}{|c|c|c|c|c|c|c|c|}
    \hline
        Exam  & Subject & Grade & Month & $\omega_1$ & $\omega_2$ & $\omega_1^s$ & $\omega_2^s$   \\
        \hline
        PBA5041F1  & Mathematics & 5th   & May   &  0.42 & 1.28 & 0.23 & 0.52 \\  
     PBA5041F2  & Mathematics & 5th   & October &  0.61 & 1.38 & 0.31 & 0.78 \\
     PBA5042F1  & Language & 5th   & May   &    0.80 & 1.61 & 0.46 & 0.73 \\ 
     PBA5042F2  & Language & 5th   & October &    0.86 & 1.51 & 0.55 & 0.95 \\
     PBA5043F1  & Science & 5th   & May   &  0.79 & 1.37 & 0.47 & 0.87 \\
     PBA5043F2  & Science & 5th   & October &   0.82 & 1.47 & 0.57 & 0.88 \\
     PBA9041F1  & Mathematics & 9th   & May   &   0.89 & 1.53 & 0.58 & 0.89 \\
     PBA9041F2  & Mathematics & 9th   & October &  1.22 & 1.53 & 0.99 & 1.07 \\
     PBA9042F1  & Language & 9th   & May   &   0.55 & 1.44 & 0.31 & 0.65 \\ 
     PBA9042F2  & Language & 9th   & October &  0.86 & 1.47 & 0.63 & 0.97 \\
     PBA9043F1  & Science & 9th   & May   &  0.78 & 1.46 & 0.59 & 0.98 \\ 
     PBA9043F2  & Science & 9th   & October &  0.78 & 1.36 & 0.63 & 1.03 \\ 
  
        \hline
        \multicolumn{8}{p{0.9\textwidth}}{Source: ICFES. Calculations: Authors.}      \\  
        \multicolumn{8}{p{0.9\textwidth}}{Number of copy-free couples accused of copying (for every 1,000 pairs) at a 99.9\% confidence level}\\
    \end{tabularx}
  \label{tab:errotipo1general2}
\end{table}

 \subsubsection{Power}\label{tipoII}
 
The following figures show the power  among the $\gamma_1$, $\gamma_1^s$, $\omega_1$ and $\omega_2^s$ indices in the Mathematics 5th grade May test. Notice that the $\omega_2^s$ index has the highest power for all levels of answer copying. This is true for all exams as shown in figures \ref{fig:comp2}-\ref{fig:comp12} in appendix \ref{power_cal}. Based on the results of the previous section and this section, we believe this justifies the use of the $\omega_2^s$ index over all otther version of the $\omega$ index and all versions of the $\gamma$ index. 

In other words, the index with the highest power among those studied, subject to the restriction of preserving the type-I error, uses a nominal response model for item answering, conditions the probability of identical answers on the answer pattern of the individual that provides answers, and calculates critical values via a normal approximation.

In the next section we apply the $\omega_2^s$ to our data and compare the prevalence of cheating across examination rooms in which different strategies to prevent cheating are used by the ICFES.

\begin{figure}[H]
	\centering
	\caption{}
		\includegraphics[height=0.35\textheight]{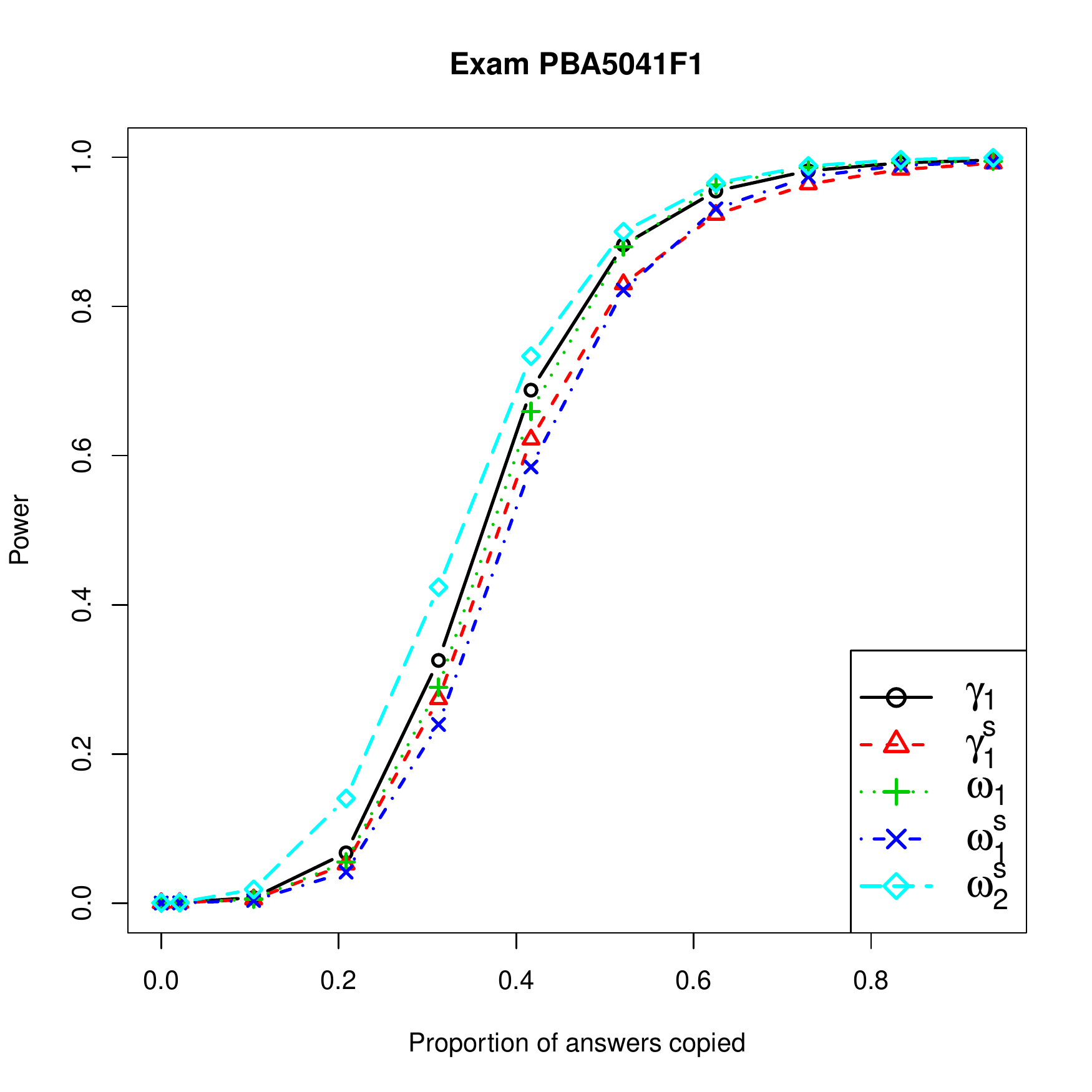}
	\caption*{Power in terms of the proportion of answers copied, for all the indices, in the Mathematics 5th grade May test.}
	\label{fig:comp1}
\end{figure}

\section{Strategies to prevent cheating}

The ICFES randomly assigns schools to three different samples that have different levels of proctoring. Most of the schools are assigned to the \textit{censal} sample in which the ICFES distributes the exams to the schools and the schools perform the proctoring. The \textit{controlada} and \textit{estadistica} samples are smaller but proctoring is done by the government itself. In the \textit{controlada} sample the proctoring is done by the central government (i.e. the ICFES) while in \textit{estadistica} sample proctoring is carried out by the regional government (Secretarias de Educaci\'{o}n). 

Table \ref{tab:1} shows some descriptive statistics of the samples. There are two things that are worth mentioning here. First, in the \textit{controlada} and \textit{estadistica} samples the ICFES had different versions of each test, so that each student was randomly assigned to one of three possible tests per subject. We only have information for the students that answered the version of the test that was used in the \textit{censal} sample. Thus one would expect the average number of students per school in the \textit{controlada} and \textit{estadistica} samples to be around one third of those in the \textit{censal} sample; however, this is not the case. Second, the October \textit{estadistica} sample has more students per school than either the \textit{controlada} or the \textit{censal} samples. These two results lead us to believe that the assignment of schools to samples was not entirely random.

\begin{figure}[H]
 \caption{Cheating across samples}
 \label{fig:comp}
\begin{center}
 \includegraphics[width=\textwidth]{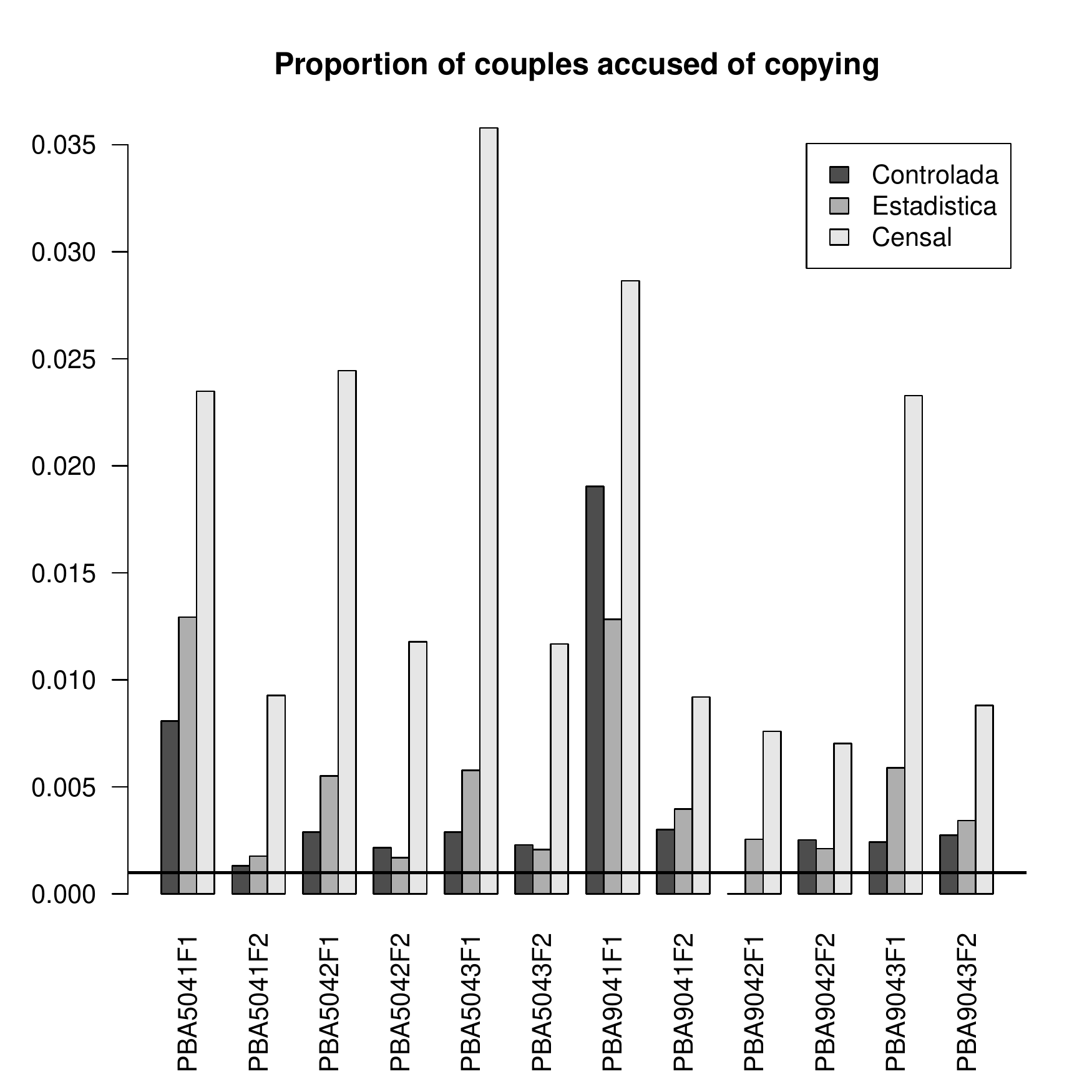}
 		\caption*{Proportion of couples accused of copying across exams  and across samples using the $\omega_2^s$ index. The horizontal black line is the theoretical type-I error rate. Source: ICFES. Calculations: Authors. }
\end{center}
 \end{figure}

We apply the $\omega_2^s$ index to the three different samples (see figure \ref{fig:comp}). In most cases the prevalence of cheating according to the index is lower for the \textit{controlada} or the \textit{estadistica} sample and highest for the \textit{censal} sample.  In most cases the \textit{controlada} sample has a lower prevalence of cheating or a similar level to the \textit{estadistica} sample, except for the May 9th grade Mathematics (PBA9041F1) test\footnote{This could be due to sampling variation given that there are only 75 schools in the \textit{controlada} sample.}.

These results can be interpreted in at least two different ways. If one remains skeptical about the index then this would serve as evidence that the index is indeed detecting cheating. Alternatively, if one believes that the index can be used as a reliable measure of cheating, these results can be interpreted as the amount of cheating that is prevented by increasing the level of proctoring. However, since the assignment of schools to samples does not seem to be random these results must be taken with caution as unobservable factors could bias the results. Additionally, it is impossible to distinguish the effect of proctoring and of having multiple versions of an exam distributed to students.

\begin{table}[H]

\caption{Characteristics of the \textit{controlada}, \textit{estadistica} and \textit{censal} samples}\label{tab:1}
\begin{center}
\begin{tabular}{lccc}
\hline
 & Controlada & Estadistica & Censal \\ 
  \hline

\multicolumn{4}{c}{5th Grade May} \\

  \hline
No. Students & 1,413 & 7,648 & 60,099 \\ 
  No. of Schools & 141 & 680 & 3,421\\ 
  Students/School & 10.02 & 11.25 & 17.57 \\ 
   & (0.88) & (0.47) & (0.46) \\ 
  \hline

\multicolumn{4}{c}{5th Grade October} \\

  \hline
No. Students & 3,830 & 26,393 & 403,624 \\ 
  No. of Schools & 958 & 654 & 31,827 \\ 
  Students/School & 4.00 & 40.36 & 12.68 \\ 
   & (0.13) & (1.36) & (0.11) \\ 

  \hline

\multicolumn{4}{c}{9th Grade May} \\
 
  \hline
No. Students & 1,150 & 6,690 & 44,577 \\ 
  No. of Schools & 75 & 351 & 1110 \\ 
  Students/School & 15.33 & 19.06 & 40.16 \\ 
   & (1.62 ) & (1.08) & (1.44) \\

  \hline

\multicolumn{4}{c}{9th Grade October} \\

  \hline
No. Students & 3,106 & 24,387 & 303,233 \\ 
  No. of Schools & 495 & 487 & 9,059 \\ 
  Students/School & 6.27 & 50.08 & 33.47 \\ 
   & (0.25) & (1.74) & (0.35) \\

  \hline
 \multicolumn{4}{l}{Source: ICFES. Calculations: Authors.  } \\ 
  \multicolumn{4}{p{0.7\textwidth}}{Note: Standard error of the mean for the number of students per school in parenthesis.} \\                           
\end{tabular}
\end{center}
\end{table}

\begin{figure}[H]
	\centering
		\caption{}
		\includegraphics[width=\textwidth]{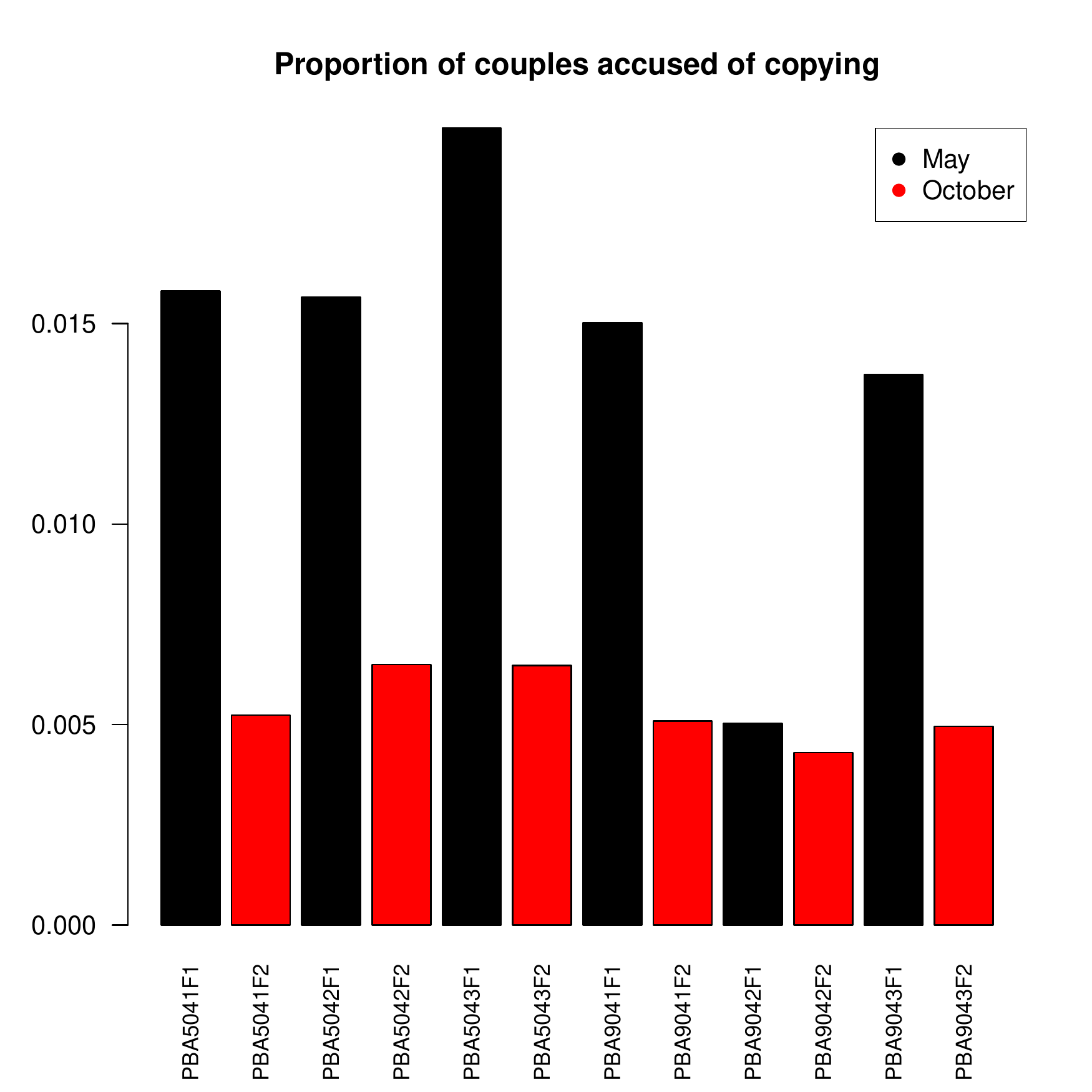}
		\caption*{Proportion of couples accused of copying across exams  in the \textit{censal} sample using the $\omega_2^s$ index. Source: ICFES. Calculations: Authors. }
	\label{fig:across}
\end{figure}

Finally, we restrict ourselves to the \textit{censal} sample. Figure \ref{fig:across} shows the proportion of couples for which the index rejects the null hypothesis. There is a clear pattern in which cheating drops dramatically between May and October. The SABER tests are administrated over three sessions, wherein students answer a different subject in each session. In May, every student took the same subject at the same time, while in October only one third of the students took the same subject in each session, thus reducing the number of students from whom one could copy in a given session. In other words, in May all students took the mathematics portion of the test at the same time. In October, while one third of the students answered the mathematics portion of the test another third answered the language portion and the final third answered the science portion. 

Note that the Language Test for 9th graders in May (PBA9042F1 test) does not follow the trend. It is also surprising to find that the levels of cheating are similar for 5th and 9th graders. These populations are different in terms of motivation, maturity and sophistication. We could not find a reasonable explanation for either of these phenomena. 

As before, these results can be interpreted in at least two different ways. They could be interpreted as evidence that the index is indeed detecting cheating. Alternatively, if one believes that the index can be used as a reliable measure of cheating, these results can be interpreted as the amount of cheating that is prevented by having different students answer different parts of the exam at different times instead of having all of them answer portion of the exam at the same time. Again this results must be taken with caution as the population of students in May and October might be different in unobservable factors which could bias the results.

%
%In short these results provide further evidence that the index is indeed measuring cheating. Additionally, if we take these results literally they provide estimates of how effective current strategies for cheating-prevention are. It is important to note that these numbers must be taken with a grain of salt as they do not come from a randomized experiment and unobserable factors could bias the results.

\section{Massive cheating}

In this section we consider a subject rarely treated in the answer-copying literature: massive cheating. Many institutions, including the ICFES, do not use answer copying-indices to formally blame an individual of copying. Rather they are interested in detecting highly suspicious examination rooms. The ICFES forces suspicious examination rooms to retake the exam under stricter surveillance conditions.

To determine whether massive cheating has occurred in an examination room, multiple hypotheses must be tested. If the significance level for a given statistical test is $\alpha_{I}$, the significance level for a multiple test ($\alpha_{MT}$) will increase exponentially as the number of hypothesis to be tested increases. In other words, $\alpha_{MT}=1-(1-\alpha_{I})^n\leq \alpha_{I} \cdot n$, assuming independence across the hypotheses. Thus we need to set $\alpha_{I}=\frac{\alpha_{MT}}{n}$ in order to assure the multiple test significance level; if this correction is made, in most cases the power of the test is severely diminished. To overcome these difficulties a line of research has developed procedures to control error rates similar to the type-I error (of the multiple hypotheses test), which can be easily applied in many cases without compromising the power of the test.

Most of these methodologies are based on Bonferroni correction that control the false positive rate (that is, the number of null hypotheses rejected incorrectly as a proportion of the number of null hypotheses rejected). We use the results from applying the $\omega_2^s$ index to every examination room. If there are $n$ students in a room then the index is applied $n\times (n-1)$ times and we adjust the p-values following the correction given by \citeA{Benjamini1995}.

%Following the notation of \citeA{Benjamini1995}, assume there are $m=n\times (n-1)$ null hypotheses to be tested (the number of couples in an examination room). The null hypothesis always states that no answer-copying has been done. Of these, $m_0$ are true and $R$ are rejected by the test. The following table summarizes the mistakes made when evaluating the $m$ hypothesis:
%
%\begin{table}[H]
%			\caption{Mistakes made when evaluating the $m$ null hypothesis}
%	\centering
%	\footnotesize
%		\begin{tabularx}{\textwidth}{|l|XXX|}		
%			    \hline
%    & Do not reject $H_0$ & Reject $H_0$ & Total \\
%    \hline
%    True null hypothesis & $U$ & $V$ & $m_0$ \\
%    False null hypothesis & $T$ & $S$ & $m-m_0$ \\
%    Total & $m-R$ & $R$ & $m$ \\
%    \hline
%		\end{tabularx}
%	\label{tab:ErroresCometidosAlProbarMHipotesisNulas}
%\end{table}
%
%Table \ref{tab:ErroresCometidosAlProbarMHipotesisNulas} can be interpreted as follows. $U$ is the number of students that did not copy and were not accused of copying by the index. $V$ is the number students that did not copy, but were accused of copying by the index. $T$ is the number of students that copied, but were not accused of copying by the index. $S$ is the number of students that copied and were accused of copying by the index. The expected false positive rate and type-I error rate of the multiple hypothesis test are $E[V/(V+S)]$ and $P(V\geq 1)$; respectively. The false positive rate and the type-I error rate are identical if all null hypotheses are true; otherwise, the false positive rate is always lower.

Suppose there are $H_1,...,H_m$ hypotheses to be tested, ordered such that their $p-values$ follow $P_1 \leq P_2 \leq ... \leq P_m$, where $P_i$ is the $p-value$ of hypothesis $H_i$. Let $k$ be the greatest integer $i$, such that:

\begin{equation}
P_i \leq \frac{i}{m} p^*.
\end{equation}

$H_i$ is then rejected for every $i \in \{1,...,k\}$. This controls for the false positive rate to a maximum of $p^*$ \cite{Benjamini1995}. The previous statement is only true if there is independence between true null hypotheses. This assumption implies that the decision to not copy is an individual one and is unrelated to the decision not to copy of other individuals. This assumption depends on the conditions under which cheating takes place. For example, if an examination room has poor supervision, it would motivate several students to copy, thus invalidating the assumption.
% If we set the false positive rate to a maximum of $p^*$ and that the null hypothesis $H_i$ (student $i$ did not copy) is rejected for all $i \in \{1,...,k\}$. There is thus massive cheating (i.e., $\frac{k}{m}\geq 60\%$) with a confidence level of $1-p^*$.

\begin{figure}[H]
	\centering
	\caption{Massive cheating per exam}
		\includegraphics[width=\textwidth]{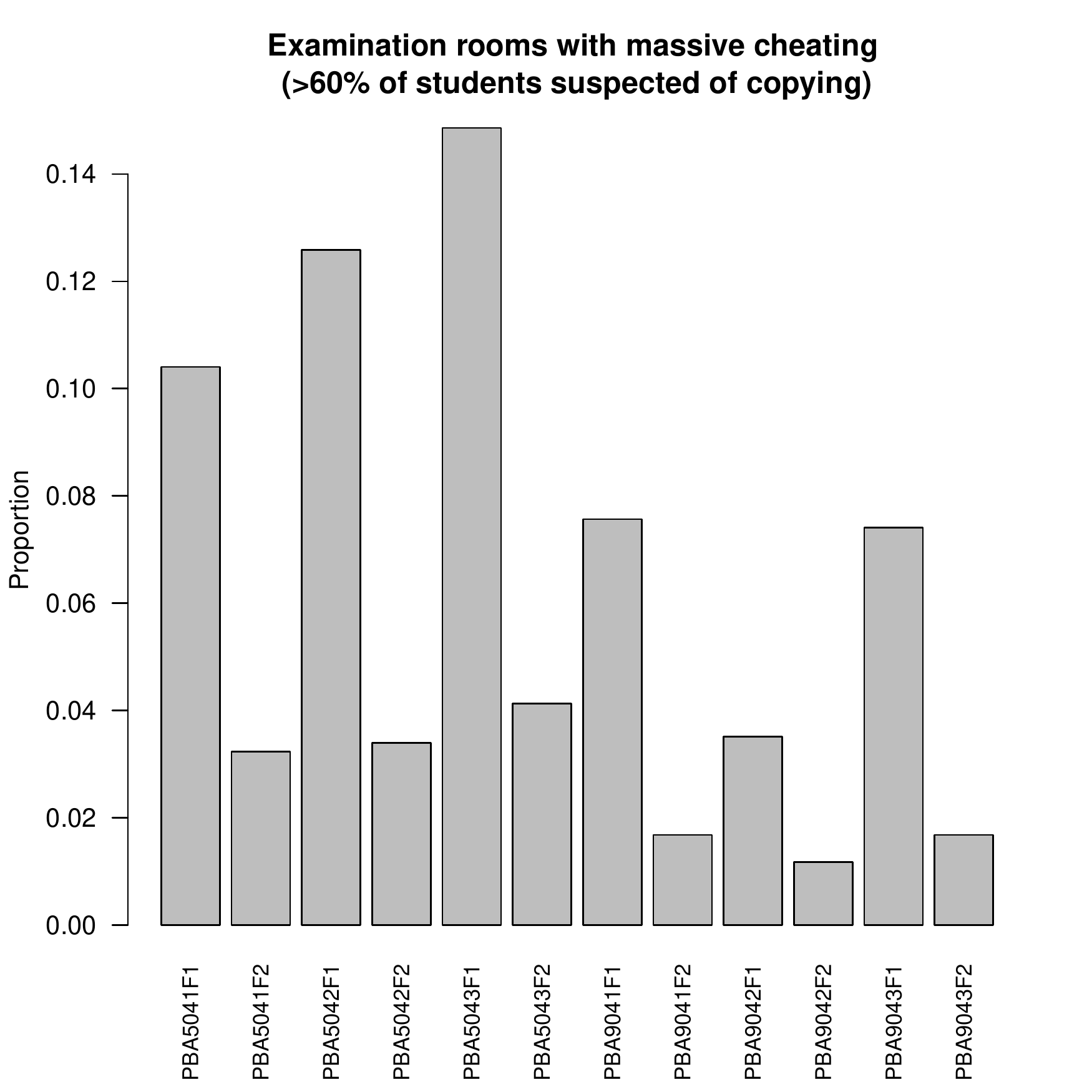}
	\caption*{Proportion of examination rooms where more than 60\% of the students are suspected of cheating after using \citeA{Benjamini1995} correction for testing multiple hypothesis and controlling the false positive rate at $p^*=0.01$. Source: ICFES. Calculations: Authors. }
	\label{fig:copia_masiva}
\end{figure}

Figure \ref{fig:copia_masiva} presents the proportion of examination rooms where more than 60\%\footnote{This is the level used by the ICFES to make students in exmaination rooms retake the test under stricter surveillance conditions.}  of the students are suspected of cheating, for $p^*=0.01\%$. As can be seen, there is a high proportion of examination rooms with massive cheating. This could be explained by the fact that the examination rooms consist of students in the same grade in a given school. Nevertheless, the fact that the proportion of massive cheating drops dramatically between May and October is reassuring, since less cheating is expected in the latter. It is also interesting that the levels of massive cheating are lower, in general, for the 9th grade\footnote{The ICFES compared our results with information they have regarding school's reputation in terms of ``honesty'', and found the two to be consistent. Unfortunately, we do not have permission to divulge this information; consequently the results of this comparison cannot be presented here.}.

%The SABER tests are administrated over three sessions, wherein students answer a different subject in each session. In May, every student took the same subject at the same time, while in October only one third of the students took the same subject in each session, thus reducing the number of students from whom one could copy in a given session. In other words, in May all students took the mathematics portion of the test at the same time. In October, while one third of the students were answering the mathematics portion of the test, another third was answering the language and the final the science portion. 

%%%%%%%%%%%%%%%%%%%%%%%
%%%%%%%%%%%%%%%%%%%%%%%%%%5
%%%%%%%%%%%%%%%%%%%%%%%%%%%%%
%%%%%%%%%%%%%%%%%%%%%%%%%%%%%%%%%
%%%%%%%%%%%%%%%%%%%%%%%%%%%%%%%%%%%%%
\section{Conclusions}

In this article we justify the use of a variety of statistical tests (known as indices) found in the literature to detect answer copying in standardized tests. We derived the uniform most powerful (UMP) test (index) using the Neyman-Pearson's Lemma under the assumption that the response distribution is known. In practice, a behavioral model for item answering must be estimated and indices vary along which model they assume.

Using data from the SABER 5th and 9th tests taken in May and October of 2009 in Colombia, we compare eight widely used indices that are based on the work of  \citeA{Frary1977,Wollack1997,Wesolowsky2000, Sotaridona2006}. We first filter out the indices that do not meet the theoretical type-I error rate in practice and then select most powerful index among them. The most powerful index, among those that respect the type-I error rate, is a conditional index that models student behavior using a nominal response model, conditions the probability of identical answers on the answer pattern of the individual that provides answers, and relies on the central limit theorem to find critical values.

Using this index we analyze 12 exams taken by 5th and 9th graders in May and October of 2009 in Colombia. We find a negative correlation between the level of proctoring and the prevalence of cheating. We also find a lower prevalence of copying in examination rooms where students answer different portions of the test at the same time compared to examination rooms where all students answer the same portion of the test at the same time. These results have at least two possible interpretations: they could be interpreted as evidence that the index is indeed detecting cheating, or, alternatively, if one believes that the index can be used as a reliable measure of cheating, these results can be interpreted as the amount of cheating that is prevented by each one of these strategies to control cheating. 

Finally, we propose a methodology for detecting massive cheating while controlling for the false positive rate using a Bonferroni correction. Institutions that use answer-copying indices should also use Bonferroni corrections to test for multiple hypothesis as this extension is straightforward. 

%Our article has four contributions. First, it justifies the use of several indices by using the Neyman-Pearson Lemma. Second, we compare the empirical type-I and type-II error rate of several indices, some of which have not previously been comparatively evaluated in the specialized literature. Third, we provide suggestive evidence on how simple controls can reduce cheating. Fourth, we provide a simple framework of a Bonferroni type procedure and an application to detect massive cheating.

We believe the results in this paper should have practical implications and lead to the use of what we call the $\omega_2^s$ over other indices and the adoption of Bonferroni corrections. Further research should be done to evaluate the effectiveness of different strategies to reduce cheating.

\bibliographystyle{apacite}
\bibliography{bib_web}

\appendix

\section{Power}\label{power_cal}

\begin{figure}[H]
	\centering
	\caption{}
		\includegraphics[height=0.35\textheight]{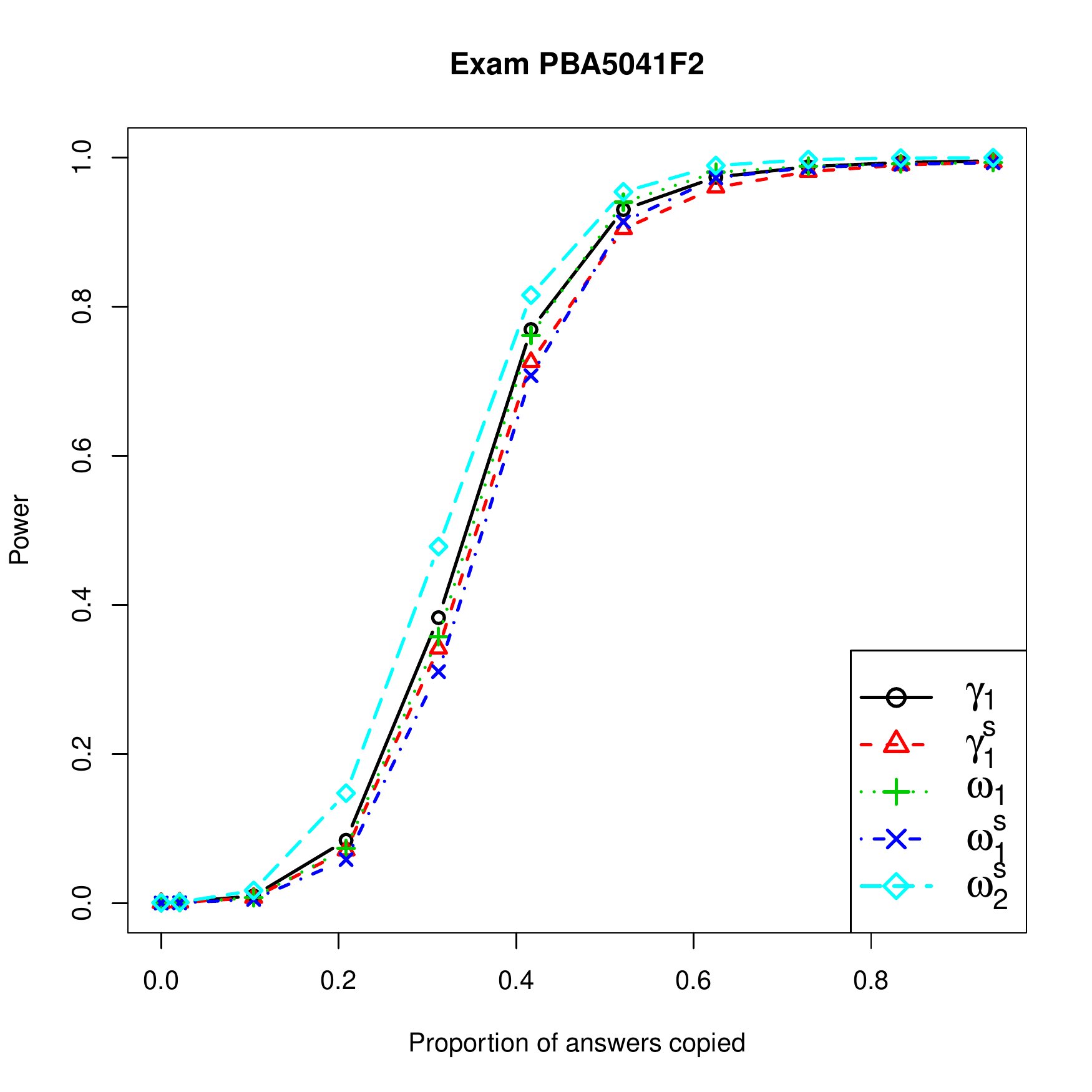}
	\caption*{Power in terms of the proportion of answers copied, for all the indices, in the mathematics 5th grade October test. Source: ICFES. Calculations: Authors.}
	\label{fig:comp2}
\end{figure}
\begin{figure}[H]
	\centering
	\caption{}
		\includegraphics[height=0.35\textheight]{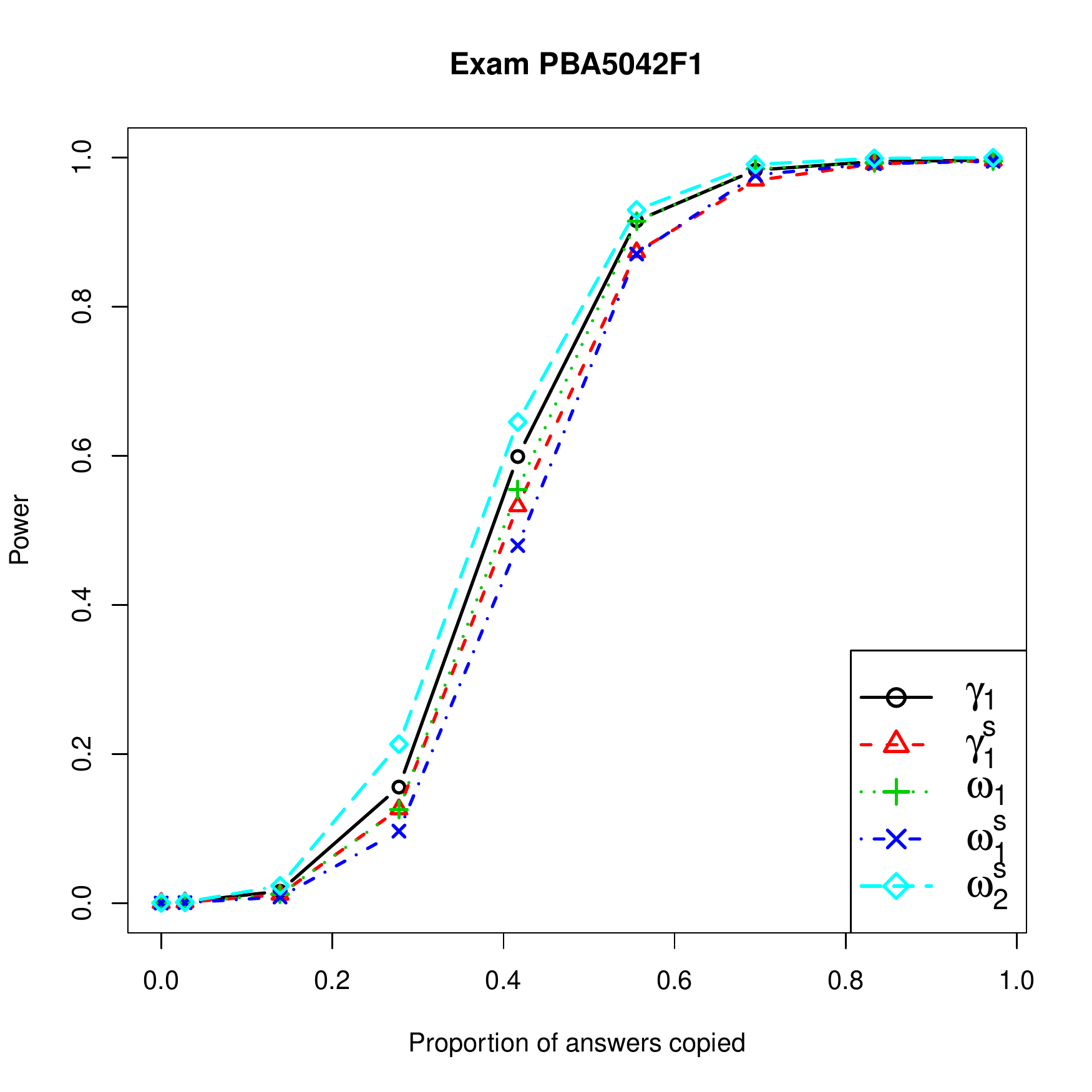}
	\caption*{Power in terms of the proportion of answers copied, for all the indices, in the language 5th grade May test. Source: ICFES. Calculations: Authors.}
	\label{fig:comp3}
\end{figure}
\begin{figure}[H]
	\centering
	\caption{}
		\includegraphics[height=0.35\textheight]{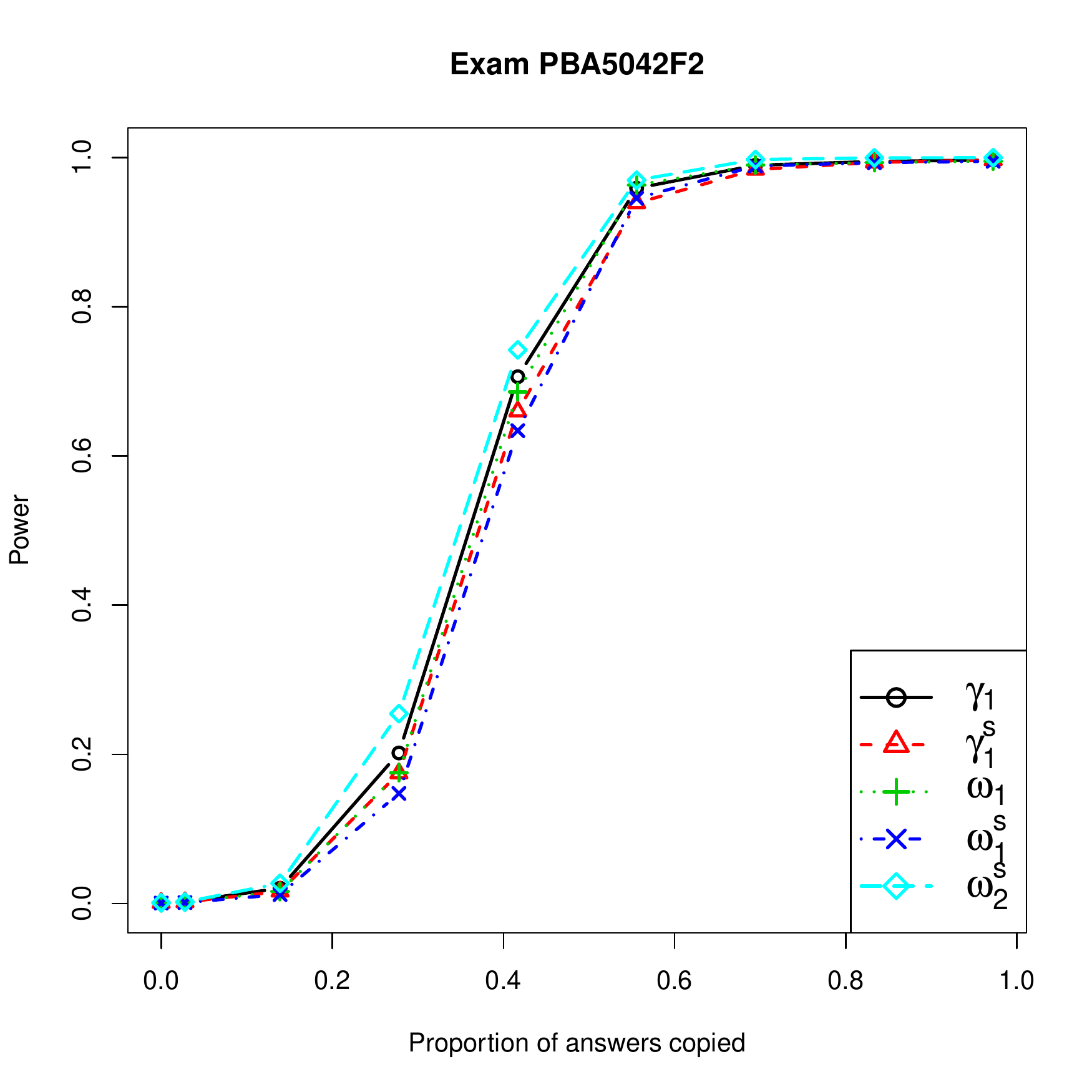}
	\caption*{Power in terms of the proportion of answers copied, for all the indices, in the language 5th grade October test. Source: ICFES. Calculations: Authors.}
	\label{fig:comp4}
\end{figure}
\begin{figure}[H]
	\centering
	\caption{}
		\includegraphics[height=0.35\textheight]{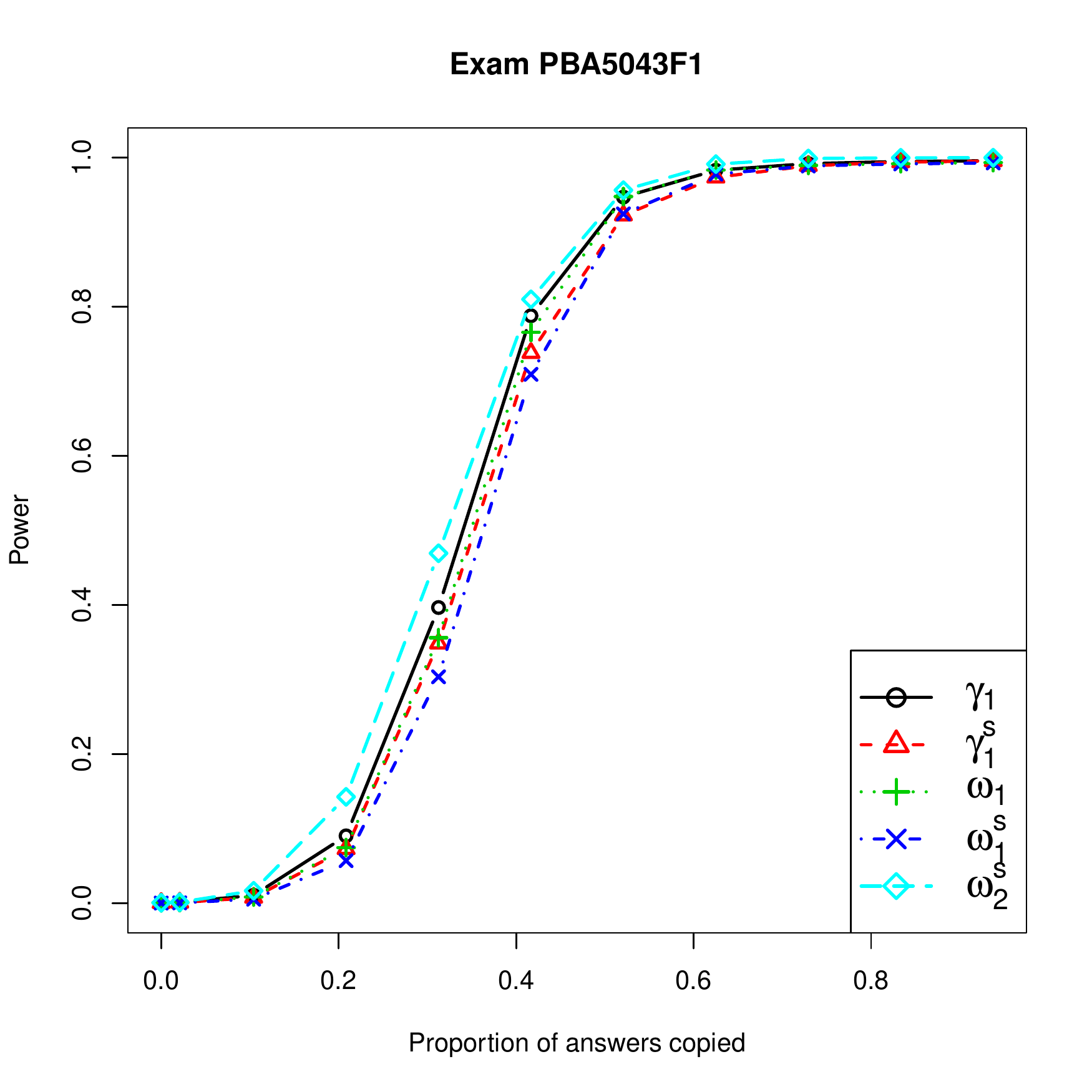}
	\caption*{Power in terms of the proportion of answers copied, for all the indices, in the science 5th grade May test. Source: ICFES. Calculations: Authors.}
	\label{fig:comp5}
\end{figure}
\begin{figure}[H]
	\centering
	\caption{}
		\includegraphics[height=0.35\textheight]{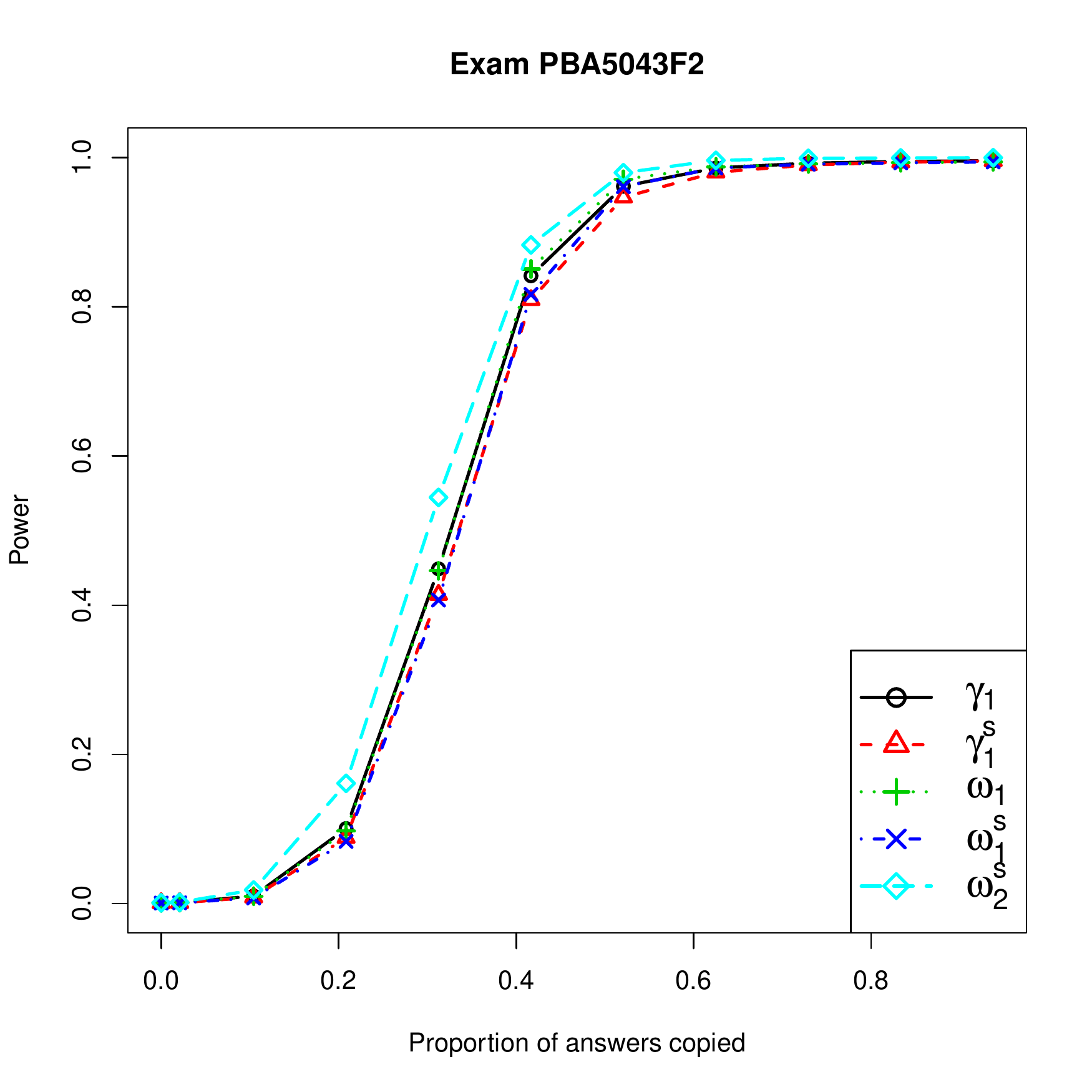}
	\caption*{Power in terms of the proportion of answers copied, for all the indices, in the science 5th grade October test. Source: ICFES. Calculations: Authors.}
	\label{fig:comp6}
\end{figure}
\begin{figure}[H]
	\centering
	\caption{}
		\includegraphics[height=0.35\textheight]{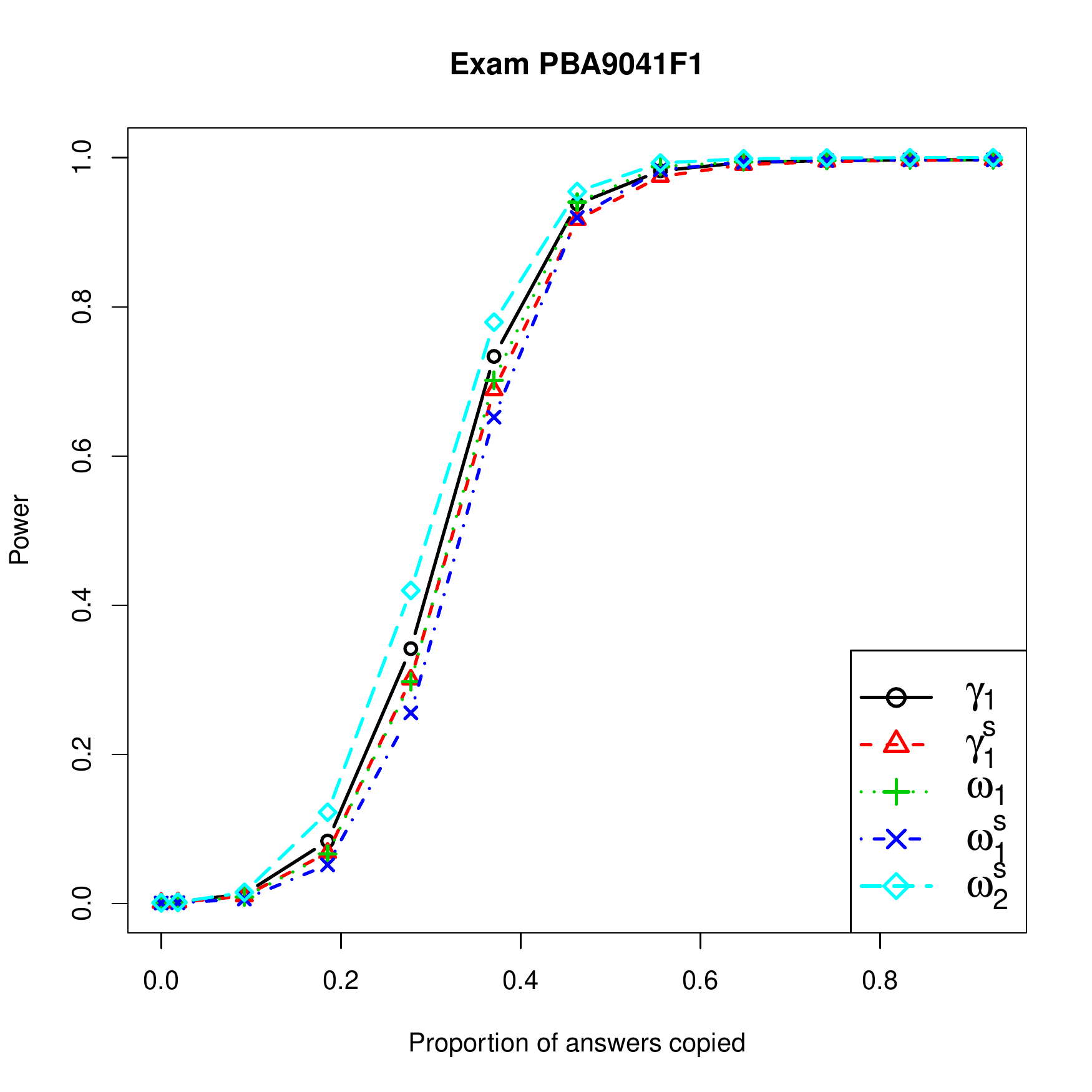}
	\caption*{Power in terms of the proportion of answers copied, for all the indices, in the mathematics 9th grade May test. Source: ICFES. Calculations: Authors.}
	\label{fig:comp7}
\end{figure}
\begin{figure}[H]
	\centering
	\caption{}
		\includegraphics[height=0.35\textheight]{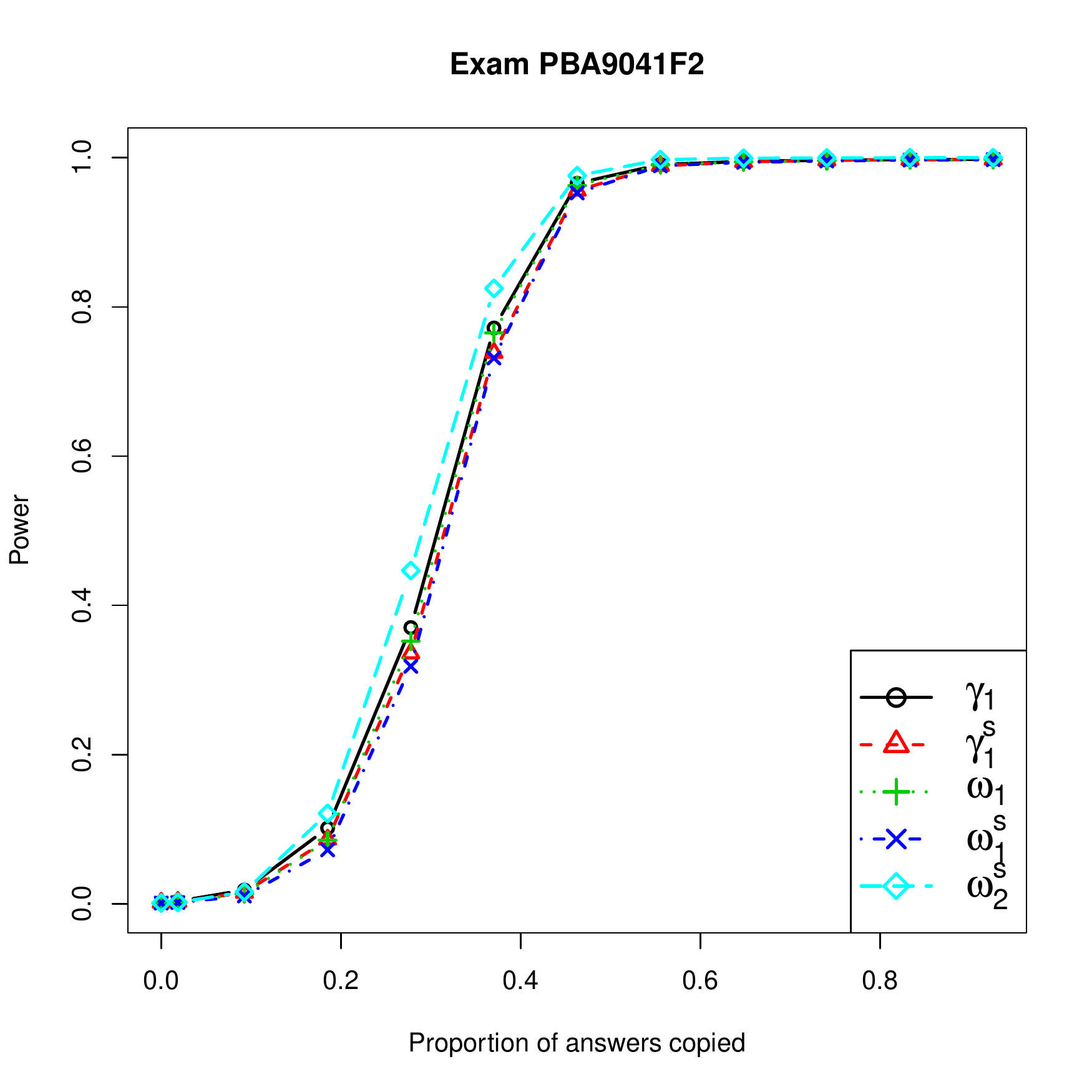}
	\caption*{Power in terms of the proportion of answers copied, for all the indices, in the mathematics 9th grade October test. Source: ICFES. Calculations: Authors.}
	\label{fig:comp8}
\end{figure}
\begin{figure}[H]
	\centering
	\caption{}
		\includegraphics[height=0.35\textheight]{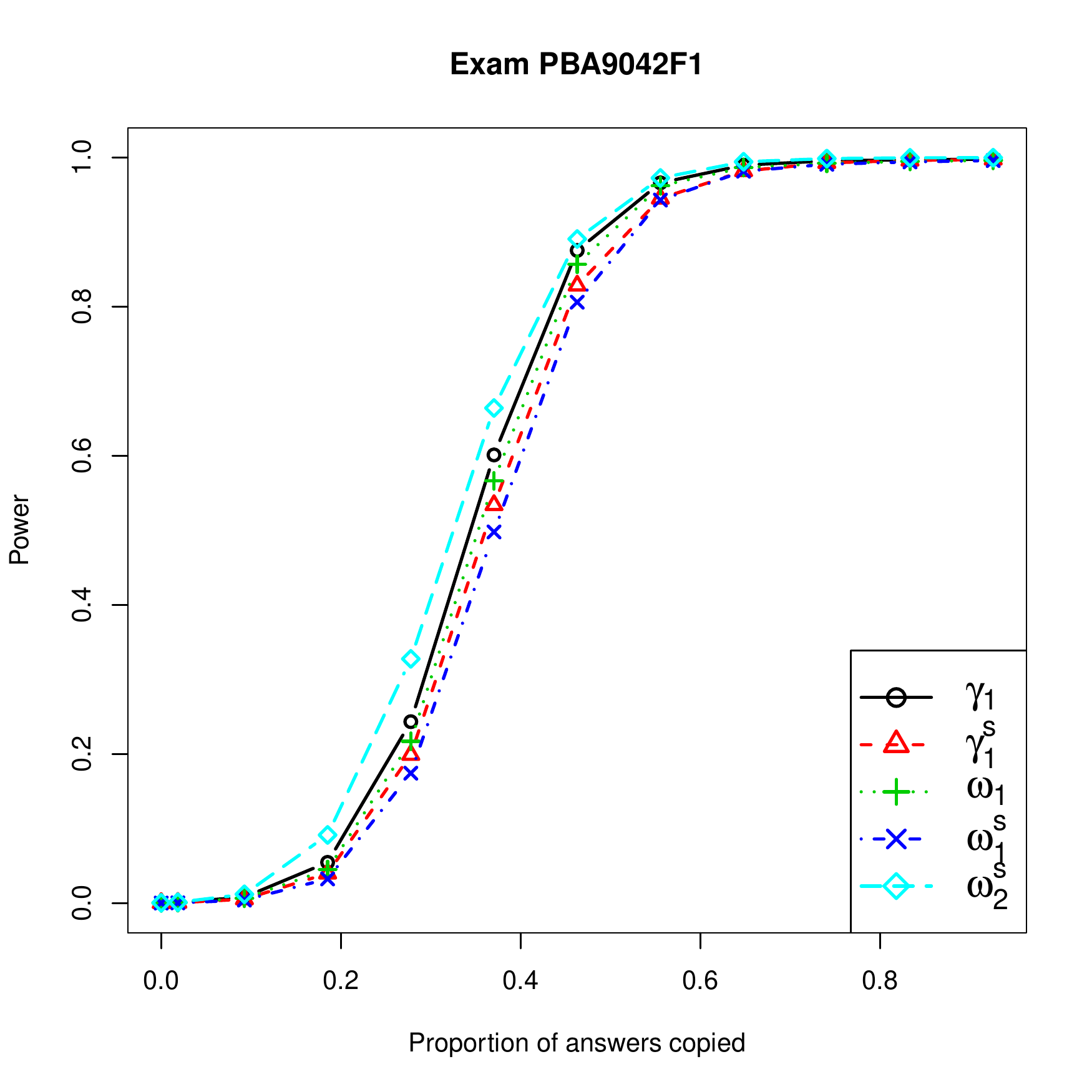}
	\caption*{Power in terms of the proportion of answers copied, for all the indices, in the language 9th grade May test. Source: ICFES. Calculations: Authors.}
	\label{fig:comp9}
\end{figure}
\begin{figure}[H]
	\centering
	\caption{}
		\includegraphics[height=0.35\textheight]{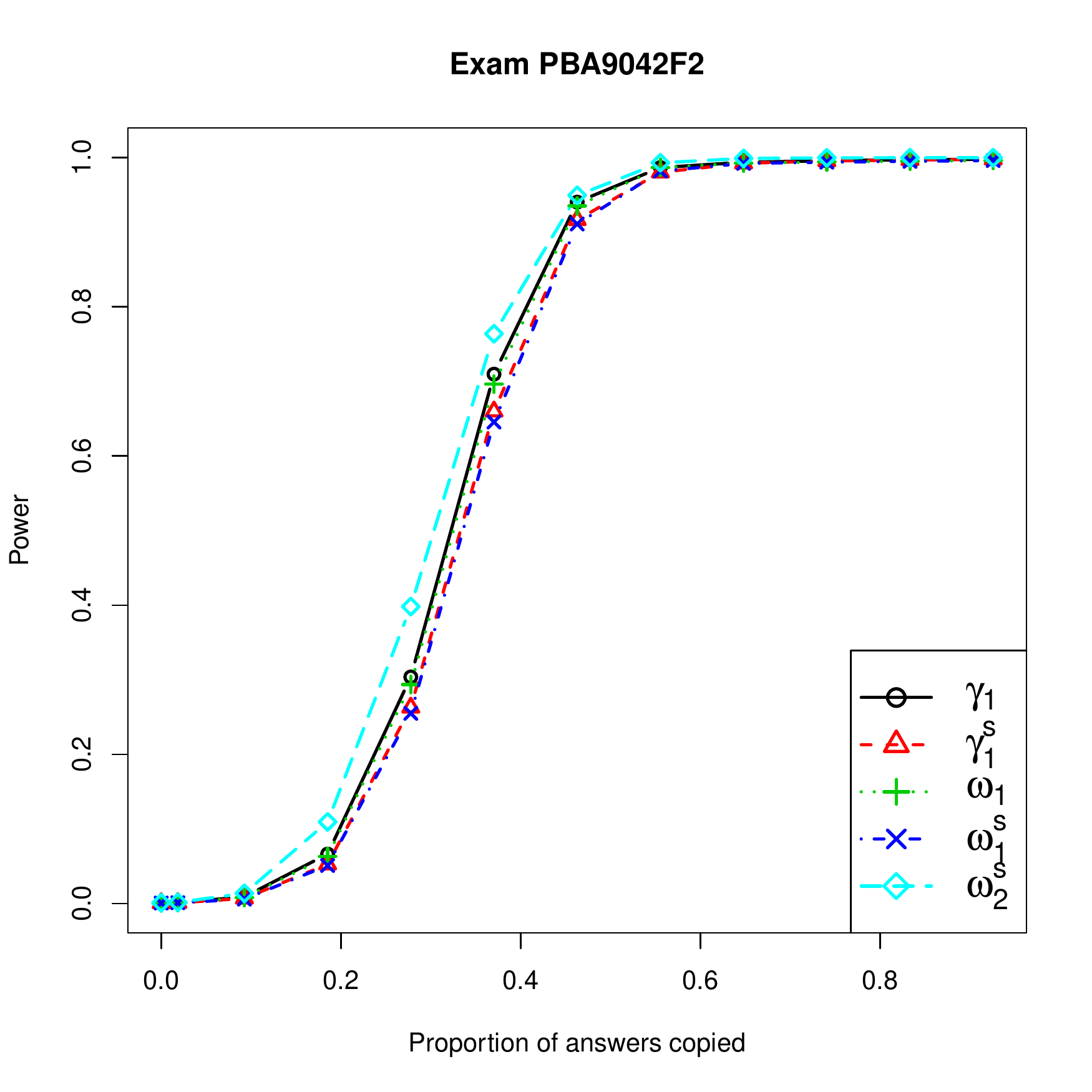}
	\caption*{Power in terms of the proportion of answers copied, for all the indices, in the language 9th grade October test. Source: ICFES. Calculations: Authors.}
	\label{fig:comp10}
\end{figure}
\begin{figure}[H]
	\centering
	\caption{}
		\includegraphics[height=0.35\textheight]{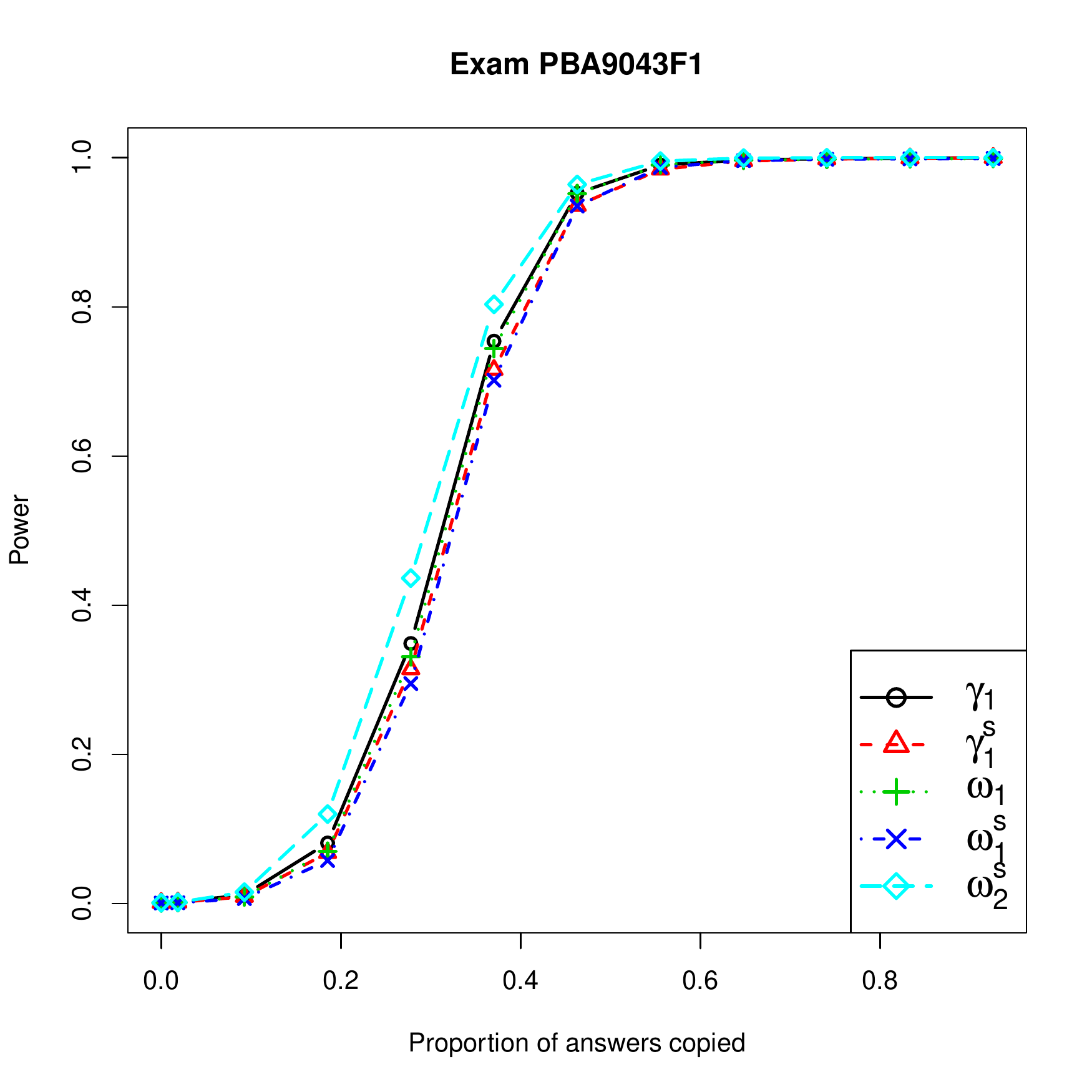}
	\caption*{Power in terms of the proportion of answers copied, for all the indices, in the science 9th grade May test. Source: ICFES. Calculations: Authors.}
	\label{fig:comp11}
\end{figure}
\begin{figure}[H]
	\centering
	\caption{}
		\includegraphics[height=0.35\textheight]{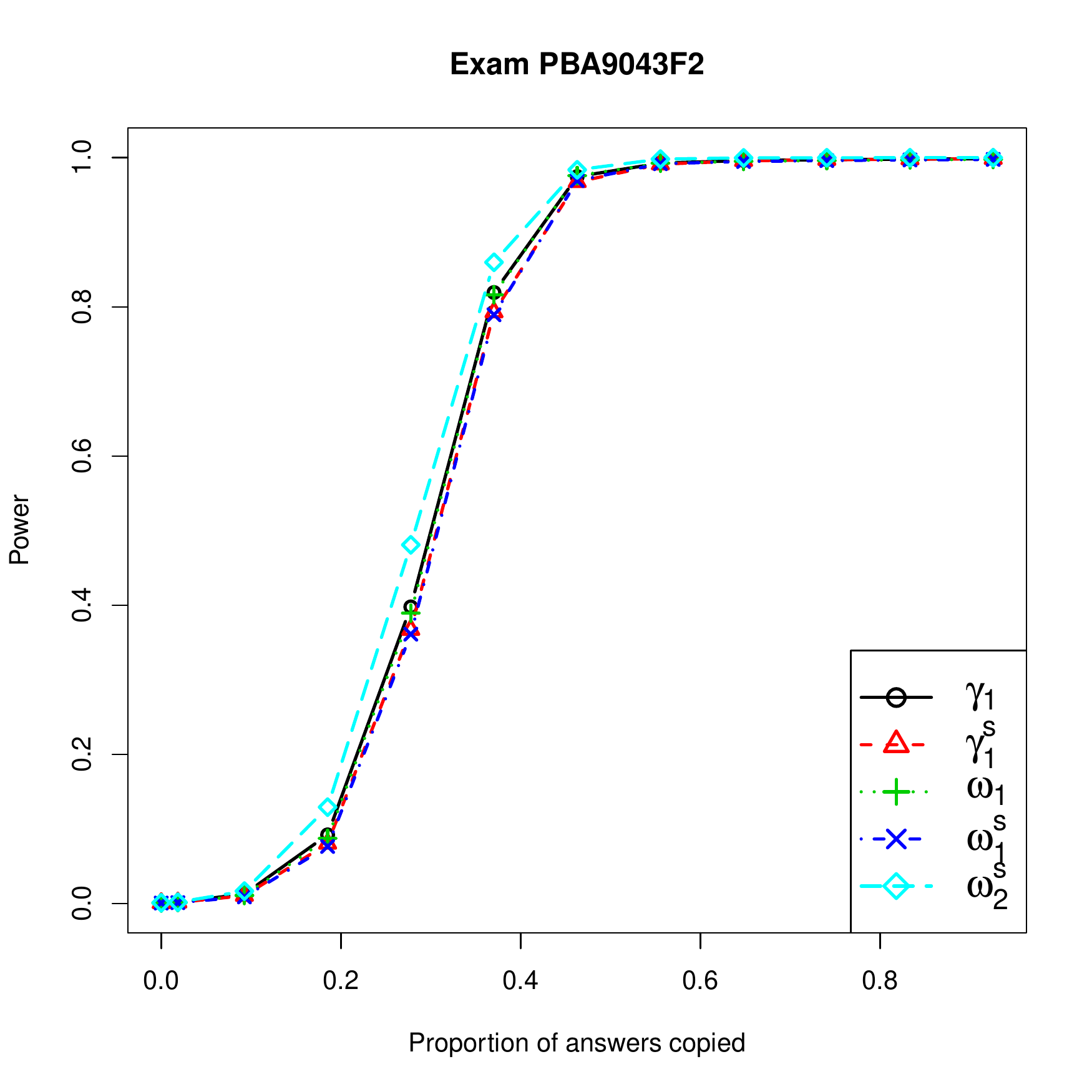}
	\caption*{Power in terms of the proportion of answers copied, for all the indices, in the science 9th grade October test. Source: ICFES. Calculations: Authors.}
	\label{fig:comp12}
\end{figure}

\end{document}